\documentclass[english]{article}
\usepackage[T1]{fontenc}
\usepackage[latin9]{inputenc}
\usepackage{color}
\usepackage{amsmath}
\usepackage{amssymb}

\makeatletter

\usepackage{amsthm}
\usepackage{fullpage}
\usepackage{makeidx}
\usepackage{amsfonts}
\usepackage{latexsym}
\numberwithin{equation}{section}

\numberwithin{figure}{section}

\newtheorem{thm}{Theorem}[section]
\newtheorem{prop}[thm]{Proposition}
\newtheorem{lemma}[thm]{Lemma}
\newtheorem{defn}[thm]{Definition}
\newtheorem{coro}[thm]{Corollary}
\newtheorem{rem}[thm]{Remark}
\newtheorem{exa}[thm]{Example}
\usepackage{babel}

\newcommand{\Z}{\mathbb{Z}}
\newcommand{\N}{\mathbb{N}}
\newcommand{\C}{\mathbb{C}}

\newcommand{\g}{\mathfrak{g}}

\newcommand{\wt}{{\rm wt}\;}

\usepackage[blocks]{authblk}
\title{On a family of vertex operator superalgebras}
\author{Haisheng Li}
\affil{ Department of Mathematical Sciences\\
 Rutgers University, Camden, NJ 08102, USA}

\author{Nina Yu}
\affil{School of Mathematical Sciences\\
 Xiamen University, Fujian 361005, CHINA}

\usepackage{babel}

\makeatother

\begin{document}
\maketitle

\abstract{This paper is to study vertex operator superalgebras which are strongly generated by their weight-$2$ and weight-$\frac{3}{2}$
homogeneous subspaces. Among the main results, it is proved that if such a vertex operator superalgebra $V$ is simple,
then $V_{(2)}$ has a canonical
commutative associative algebra structure equipped with a non-degenerate symmetric associative bilinear form and
$V_{(\frac{3}{2})}$ is naturally a $V_{(2)}$-module equipped with a $V_{(2)}$-valued symmetric bilinear form and
a non-degenerate ($\C$-valued) symmetric bilinear form, satisfying a set of conditions.
On the other hand, assume that $A$ is any commutative associative algebra equipped with a non-degenerate symmetric associative bilinear form and assume that $U$ is an $A$-module equipped with a symmetric $A$-valued bilinear form
and a non-degenerate ($\C$-valued) symmetric  bilinear form, satisfying the corresponding conditions.
Then we construct a Lie superalgebra $\mathcal{L}(A,U)$ and a simple vertex operator superalgebra
$L_{\mathcal{L}(A,U)}(\ell,0)$ for every nonzero number $\ell$
such that $L_{\mathcal{L}(A,U)}(\ell,0)_{(2)}=A$ and $L_{\mathcal{L}(A,U)}(\ell,0)_{(\frac{3}{2})}=U$.  }

\section{Introduction}

The main purpose of this paper is to determine the structure of a particular family of vertex operator superalgebras.
Note that vertex operator algebras are a special family of vertex algebras.
From definition, a vertex operator algebra $V$ is $\Z$-graded by the conformal weights: $V=\oplus_{n\in \Z}V_{(n)}$, where
$V_{(n)}$ coincides with the $L(0)$-eigenspace with eigenvalue (namely conformal weight) $n$.
The following is a basic property:
\begin{align}\label{Z-graded-va-def}
v_mV_{(k)}\subset V_{(k+n-m-1)}\quad \text{for }v\in V_{(n)},\ n, m,k\in \Z.
\end{align}
This leads to a notion of $\Z$-graded vertex algebra, where a $\Z$-graded vertex algebra is a vertex algebra
$V$ equipped with a $\Z$-grading $V=\oplus_{n\in \Z}V_{(n)}$ such that (\ref{Z-graded-va-def}) holds.
By definition (at least for this paper), a vertex operator superalgebra is
$\frac{1}{2}\Z$-graded by the conformal weights, where the even (odd) part equals
the sum of homogenous subspaces of (half) integer conformal weights.
The notion of $\frac{1}{2}\Z$-graded vertex superalgebra is defined in the obvious way.

Among the better known examples, affine vertex algebras (cf. \cite{FZ}, \cite{Lian}) can be characterized as $\Z$-graded vertex algebras
$V=\oplus_{n\in \Z}V_{(n)}$ which satisfy the conditions that $V_{(n)}=0$ for all negative integers $n$,  $V_{(0)}=\C {\bf 1}$ and
that $V_{(1)}$ generates $V$ as a vertex algebra. Note that a vertex algebra $V$ is generated by a subset $T$ if
$V$ is linearly spanned by the vacuum vector ${\bf 1}$ and vectors
$u^{(1)}_{m_1}u^{(2)}_{m_2}\cdots u^{(r)}_{m_r}{\bf 1}$ for $r\ge 1,\ u^{(i)}\in T,\ m_i\in \Z$.
A vertex (super)algebra $V$ is said to be {\em strongly generated} by a subset $T$ if $V$ is linearly spanned by
the vacuum vector ${\bf 1}$ and vectors
$u^{(1)}_{-n_1}u^{(2)}_{-n_2}\cdots u^{(r)}_{-n_r}{\bf 1}$ for $r\ge 1,\ u^{(i)}\in T,\ n_i\ge 1$
(see \cite{Kac}).
In terms of this notion, affine vertex algebras can be characterized as $\Z$-graded vertex algebras
which are strongly generated by their degree-one subspaces.

In \cite{Lam}, Lam studied (essentially) vertex operator algebras which are strongly generated by their weight-two subspaces.
 Among the main results, Lam proved that for such a vertex operator algebra $V=\oplus_{n\in \Z}V_{(n)}$,
 $V_{(2)}$ has a canonical commutative associative algebra structure with a symmetric and associative bilinear form.
From another direction, Primc (see \cite{Pr}) associated $\Z$-graded vertex algebras to Novikov algebras,
where a (left) Novikov algebra is a non-associative algebra $\mathcal{A}$ such that
 \begin{align*}
(ab)c-a(bc)=(ba)c-b(ac),\quad\quad (ab)c=(ac)b\quad \text{ for }a,b,c\in \mathcal{A}.
 \end{align*}
(Notice that a Novikov algebra with a left identity is the same as a commutative associative algebra with identity.)
For any Novikov algebra $\mathcal{A}$ with a symmetric bilinear form $\langle\cdot,\cdot\rangle$ such that
\begin{align}\label{In-0}
\langle ab,c\rangle=\langle a,bc\rangle,\ \ \ \langle ab,c\rangle=\langle ba,c\rangle\ \ \mbox{ for }a,b,c\in \mathcal{A},
\end{align}
Primc constructed a Lie algebra $\widetilde{L}(\mathcal{A}):=\mathcal{A}\otimes \C[t,t^{-1}]+\C K$ and
then a vertex algebra $V_{\widetilde{L}(\mathcal{A})}(\ell,0)$ for each complex number $\ell$.
It was shown in \cite{BLP} that these vertex algebras are exactly all the $\Z$-graded vertex algebras
which are strongly generated by their degree-two subspaces.

 In this current paper,  we study vertex operator superalgebras $V=\oplus_{q\in\frac{1}{2}\mathbb{Z}}V_{(q)}$
which are  strongly generated by their subspace $V_{(\frac{3}{2})}+V_{(2)}$.
Among the main results, we prove that if $V$ is a simple vertex operator superalgebra
which is strongly generated by its subspace $V_{(\frac{3}{2})}+V_{(2)}$,
then $V_{(2)}$ is naturally a commutative associative algebra with a non-degenerate
symmetric associative bilinear form and $V_{(\frac{3}{2})}$ is a $V_{(2)}$ module with
a non-degenerate symmetric  bilinear form and a $V_{(2)}$-valued symmetric bilinear form,
satisfying certain compatibility conditions.
On the other hand, from any commutative associative algebra $A$ with a non-degenerate
symmetric associative bilinear form and an $A$-module $U$ with
a non-degenerate symmetric  bilinear form and an $A$-valued symmetric bilinear form,
satisfying the corresponding compatibility conditions,
we construct a simple vertex operator superalgebra $V$ which is strongly generated by
its subspace $V_{(\frac{3}{2})}+V_{(2)}$, such that
$V_{(2)}\simeq A$ as an algebra and $V_{(\frac{3}{2})}\simeq U$ as an $A$-module.

Now, we continue the introduction to give some detailed information.
Let $V=\oplus_{q\in\frac{1}{2}\mathbb{Z}}V_{(q)}$ be a general vertex operator superalgebra which is strongly generated by
its subspace $V_{(\frac{3}{2})}+V_{(2)}$.
We first show that $V_{(q)}=0$ for $q\in -\frac{1}{2}\Z_{+}$, $V_{(0)}=\C {\bf 1}$, $V_{(\frac{1}{2})}=0=V_{(1)}$.
 Then we show that $V_{(2)}$ is a commutative non-associative algebra with $ab=a_{1}b$ for $a,b\in V_{(2)}$
 and the bilinear form $\langle\cdot,\cdot\rangle$ on $V_{(2)}$ defined by $a_{3}b=\langle a,b\rangle\boldsymbol{1}$
is symmetric and associative. Furthermore, we define  an action of $V_{(2)}$ on $V_{(\frac{3}{2})}$ by
$$a\cdot v=\frac{4}{3}a_{1}v\quad \text{ for }a\in V_{(2)},\ v\in V_{(\frac{3}{2})}.$$
Even though $V_{(2)}$ is only a non-associative algebra in general, its action on
$V_{(\frac{3}{2})}$ is proved to be associative in the sense that $a(bv)=(ab)v$ for $a,b\in V_{(2)},\ v\in  V_{(\frac{3}{2})}$.
On the other hand, we equip $V_{(\frac{3}{2})}$ with a $\mathbb{C}$-valued bilinear form
$\langle\cdot,\cdot\rangle$ and a $V_{(2)}$-valued bilinear form
$(\cdot,\cdot)$ on $V_{(\frac{3}{2})}$, which are defined by
$$u_{2}v=\langle u,v\rangle\boldsymbol{1}\   \text{ and }\  (u,v)=u_{0}v \  \text{  for }u,v\in V_{(\frac{3}{2})}.$$
We show that the defined bilinear forms on $V_{(2)}$ and $V_{(\frac{3}{2})}$ satisfy a set of compatibility conditions
(see (\ref{main-conditions}) for details).
On the other hand, we show that for any such vertex operator superalgebra $V$,
there exists a unique symmetric invariant bilinear form $\langle\cdot,\cdot\rangle$ with $\langle {\bf 1}, {\bf 1}\rangle=1$
and the kernel of $\langle\cdot,\cdot\rangle$ coincides with the maximal ideal of $V$.
Then by making use of this invariant bilinear form we prove that if such a vertex operator superalgebra $V$ is simple,
then $V_{(2)}$ is a commutative associative algebra and
the ($\C$-valued) bilinear forms on $V_{(2)}$ and $V_{(\frac{3}{2})}$ are non-degenerate.

In the other direction, we show that the aforementioned structures on $V_{(2)}$ and $V_{(\frac{3}{2})}$
are also sufficient to determine such a simple vertex operator superalgebra $V$.
More specifically, from any commutative associative algebra $A$, equipped with a symmetric bilinear form
$\langle\cdot,\cdot\rangle$, and an $A$-module $U$, equipped with a $\mathbb{C}$-valued bilinear form
$\langle\cdot,\cdot\rangle$ and an $A$-valued bilinear form $(\cdot,\cdot)$ on $U$,
satisfying the corresponding conditions,
we construct a vertex operator superalgebra $V$ with $V_{(2)}=A$ and $V_{(\frac{3}{2})}=U$ (under a canonical identification).
Assuming that the ($\C$-valued) bilinear forms on $A$ and $U$ are non-degenerate,
we furthermore obtain a simple vertex operator superalgebra
{\em still} with $A$ and $U$ as its weight-$2$ and weight-$\frac{3}{2}$ homogeneous subspaces.
To achieve this, we first construct a Lie superalgebra  $\mathcal{L}(A,U)$ with
$$ \mathcal{L}(A,U)=  A\otimes \C[t,t^{-1}]\oplus U\otimes \C[t,t^{-1}]\oplus \C K$$
as a vector space, where $K$ is a nonzero central element (see Section 4 for details). Then just as  affine Lie algebras,
 for every complex number $\ell$ we construct a $\frac{1}{2}\Z$-graded vertex superalgebra $V_{\mathcal{L}(A,U)}(\ell,0)$,
whose underlying space is a particular generalized Verma module for the Lie superalgebra $\mathcal{L}(A,U)$
with $K$ acting as scalar $\ell$. We furthermore obtain the desired simple vertex superalgebras as quotients
of $V_{\mathcal{L}(A,U)}(\ell,0)$ for nonzero $\ell$.
By using generalized Verma modules for the Lie superalgebra $\mathcal{L}(A,U)$,
we also determine all irreducible $\frac{1}{2}\N$-graded $V_{\mathcal{L}(A,U)}(\ell,0)$-modules,
which are parameterized by $\lambda\in A^{*}$.

This paper is organized as follows. We first review basic notions and facts about
vertex operator superalgebras and their modules in Section 2.
In Section 3, we determine vertex operator superalgebras $V$
which are strongly generated by their subspace $V_{(\frac{3}{2})}+V_{(2)}$. More specifically,
we show that $V_{(2)}$ has a canonical commutative associative algebra structure
with a non-degenerate symmetric associative bilinear form while
$V_{(\frac{3}{2})}$ is a $V_{(2)}$-module structure, and we give a set of compatibility conditions among these datum.
In Section 4,  from any commutative associative algebra
$A$  and an $A$-module $U$ that satisfy certain conditions, we construct a Lie superalgebra  $\mathcal{L}(A,U)$.
In Section 5, we associate vertex (operator) superalgebras
$V_{\mathcal{L}(A,U)}(\ell,0)$ to the Lie superalgebra  $\mathcal{L}(A,U)$ obtained in Section 4 and obtain simple vertex (operator) superalgebras $L_{\mathcal{L}(A,U)}(\ell,0)$.
We also determine all irreducible $\frac{1}{2}\mathbb{N}$-graded $V_{\mathcal{L}(A,U)}(\ell,0)$-modules.

In this paper, we work on the field $\mathbb{C}$ of complex numbers; all vector spaces are assumed to be over $\C$,
and we use $\mathbb{N}$ for the set of nonnegative integers and $\mathbb{Z}_{+}$
for the set of positive integers.

\section{Basics}

We start by reviewing some classical notions. A {\em superalgebra}
is a non-associative algebra $A$ with a $\mathbb{Z}_{2}$-grading
$A=A_{\bar{0}}\oplus A_{\bar{1}}$ such that $A_{\alpha}\cdot A_{\beta}\subset A_{\alpha+\beta}$
for $\alpha,\beta\in\Z_{2}$. If a superalgebra $A$ has an identity
$1$, then $1\in A_{\bar{0}}$. For a superalgebra $A$, a (homogeneous)
element of $A_{\bar{0}}$ is said to be {\em even} while an element
of $A_{\bar{1}}$ is said to be {\em odd}. For a nonzero homogeneous
element $a$ of $A$, define $\left|a\right|=0$ if $a$ is even and
$\left|a\right|=1$ if $a$ is odd. When there are other gradings,
for clarity we also use the term ``$\Z_{2}$-homogeneous vector.''

\begin{defn}
{\em A {\em Lie superalgebra} is a superalgebra $\g=\g_{\bar{0}}\oplus\g_{\bar{1}}$ with multiplication
$[\cdot,\cdot]$, satisfying the following conditions:
\begin{eqnarray}
 &  & [x,y]=-(-1)^{|x||y|}[y,x],\\
 &  & [x,[y,z]]-(-1)^{|x||y|}[y,[x,z]]=[[x,y],z]
\end{eqnarray}
for $x,y,z\in\g$ with $x,y\in \g_{\bar{0}}\cup \g_{\bar{1}}$ (homogeneous).}
\end{defn}

Next, we recall the definition of a vertex operator superalgebra (see
\cite{FFR}, \cite{Li}, \cite{Xu}, \cite{Kac}).

\begin{defn}
{\em A \emph{vertex operator superalgebra} is a quadruple
$(V,Y,\boldsymbol{1},\omega)$, which consists of a $\frac{1}{2}\Z$-graded
vector space $V=\oplus_{q\in\frac{1}{2}\Z}V_{(q)}$, viewed as a
superspace $V=V^{\bar{0}}\oplus V^{\bar{1}}$ with
\begin{align}
V^{\bar{0}}=\oplus_{n\in\Z}V_{(n)},\quad V^{\bar{1}}=\oplus_{n\in\Z}V_{(n+\frac{1}{2})},
\end{align}
a linear map
\begin{eqnarray*}
Y(\cdot,z):\  V & \to & (\text{End}V)[[z,z^{-1}]]\\
v & \mapsto &  Y(v,z)=\sum_{n\in\mathbb{Z}}v_{n}z^{-n-1}\ \ (\text{where }v_{n}\in\text{End}V),
\end{eqnarray*}
and two distinguished even vectors $\boldsymbol{1},$ called the {\em
vacuum vector,} and $\omega$, called the {\em conformal or Virasoro}
vector, satisfying the following conditions:

(1) For $u\in V^{\alpha}$, $v\in V^{\beta}$ with $\alpha,\beta\in\Z_{2}$,
$u_{n}v\in V^{\alpha+\beta}$ for all $n\in\mathbb{Z}$, and $u_{n}v=0$
for $n$ sufficiently large.

(2) For ($\Z_{2}$-homogeneous) vectors $u,v\in V^{\bar 0}\cup V^{\bar 1}$, the following
identity holds on $V$:
\begin{eqnarray}\label{Jacobi-delta}
&&z_0^{-1}\delta\left(\frac{z_1-z_2}{z_0}\right)Y(u,z_1)Y(v,z_2)
-(-1)^{|u||v|}z_0^{-1}\delta\left(\frac{z_2-z_1}{-z_0}\right)Y(v,z_2)Y(u,z_1)\nonumber\\
&&\quad\hspace{2cm} =z_2^{-1}\delta\left(\frac{z_1-z_0}{z_2}\right)Y(Y(u,z_0)v,z_2)
\end{eqnarray}
(the {\em Jacobi identity}).

(3) Vector ${\bf 1}$ satisfies the following vacuum and the creation properties
\[
\boldsymbol{1}_{m}v=\delta_{m,-1}v,\ \ v_{-1}\boldsymbol{1}=v\ \ \text{and }\ v_{n}\boldsymbol{1}=0\ \ \text{ for }v\in V,\ m\in\mathbb{Z},\ n\in\mathbb{N}.
\]

(4) The assignment $L(n)=\omega_{n+1}$ for $n\in\Z$ gives a representation
of the Virasoro algebra on $V$ of a central charge $c_{V}\in\mathbb{C},$ i.e.,
\begin{eqnarray}
[L(m),L(n)]=(m-n)L(m+n)+\frac{m^{3}-m}{12}\delta_{m+n,0}c_{V}
\end{eqnarray}
for $m,n\in\mathbb{Z}$.

(5) For any $v\in V$,
\begin{align*}
Y(L(-1)v,z)=\frac{d}{dz}Y(v,z),
\end{align*}
and for $q\in\frac{1}{2}\Z$,
\begin{eqnarray}
V_{(q)}=\{v\in V\ |\ L(0)v=qv\}.
\end{eqnarray}
It is also assumed that $\dim V_{(q)}<\infty$ for $q\in\frac{1}{2}\Z$
and $V_{(q)}=0$ for $q$ sufficiently negative. For $v\in V_{(q)}$ with $q\in \frac{1}{2}\Z$, we say $v$ is conformally
or $\frac{1}{2}\Z$-homogeneous of {\em (conformal) weight} $q$
and write $\text{wt}v=q$. }
\end{defn}

We also define a notion of {\em conformal vertex superalgebra} by using all the axioms above except
the last assumption on grading.

Let $V$ be a conformal vertex  superalgebra.
In terms of components, Jacobi identity (\ref{Jacobi-delta}) for $\Z_{2}$-homogeneous vectors $u,v\in V$ can be written as
\begin{eqnarray}
\sum_{i=0}^{\infty}(-1)^{i}\binom{m}{i}\left(u_{m+k-i}v_{n+i}-(-1)^{|u||v|}(-1)^{m}v_{m+n-i}u_{k+i}\right)
=\sum_{i=0}^{\infty}\binom{k}{i}\left(u_{m+i}v\right)_{n+k-i}
\end{eqnarray}
for $m,n,k\in\mathbb{Z}$.  As immediate consequences, we have
\begin{equation}\label{commutator formula for all cases}
u_{m}v_{n}-\left(-1\right)^{|u||v|}v_{n}u_{m}=\sum_{i=0}^{\infty}\binom{m}{i}\left(u_{i}v\right)_{m+n-i}
\end{equation}
(the {\em Borcherds commutator formula}) and
\begin{equation}
\left(u_{m}v\right)_{n}=\sum_{i=0}^{\infty}(-1)^{i}\binom{m}{i}\left(u_{m-i}v_{n+i}-(-1)^{|u||v|}\left(-1\right)^{m}v_{m+n-i}u_{i}\right)\label{asso commutator}
\end{equation}
(the {\em Borcherds iterate formula}) for $m,n\in\Z$.
We also have
\begin{eqnarray}\label{skew-symmetry}
Y(u,z)v=(-1)^{|u||v|}e^{zL(-1)}Y(v,-z)u
\end{eqnarray}
(the {\em skew-symmetry}), which amounts to
\begin{equation}
u_{n}v=(-1)^{|u||v|}(-1)^{n+1}\left(v_{n}u+\sum_{i\ge1}\frac{1}{i!}(-1)^{i}L(-1)^{i}v_{n+i}u\right)\label{skew-symmetry-component}
\end{equation}
for $n\in\mathbb{Z}$.

Recalling that
\begin{eqnarray}
L(n)=\omega_{n+1}\ \ \ \ \text{ for }n\in\Z,
\end{eqnarray}
we have
\begin{eqnarray}
[L(-1),v_{n}]=(L(-1)v)_{n}=-nv_{n-1}\ \ \ \ \text{ for }v\in V,\ n\in\Z,
\end{eqnarray}
\begin{eqnarray}
u_{n}V_{(q)}\subset V_{(p+q-n-1)}
\end{eqnarray}
for $u\in V_{(p)},\ p,q\in\frac{1}{2}\Z,\ n\in\Z$.

A nonzero vector $v$ of $V$ is called a {\em quasi primary vector
of weight $h\in\frac{1}{2}\Z$} (see \cite{FHL}) if $v\in V_{(h)}$
and $L(1)v=0$. If $v$ is a quasi primary vector of weight $h$, then
\begin{eqnarray}
[L(1),v_{n}]=(2h-n-2)v_{n+1}\quad\text{ for }n\in\Z.\label{quasi-primary}
\end{eqnarray}

\begin{defn} {\em Let $V$ be a conformal vertex superalgebra. A
{\em weak $V$-module} is a vector space $M$ equipped with a linear
map
\begin{align*}
Y_{M}(\cdot,z):\  & V\to({\rm End}\;M)[[z,z^{-1}]]\\
 & v\mapsto Y_{M}(v,z)=\sum_{n\in\mathbb{Z}}v_{n}z^{-n-1},
\end{align*}
satisfying the following conditions:

(1) For any $v\in V,\ w\in M$, $v_{n}w=0$ for $n$ sufficiently positive.

(2) $Y_{M}({\bf 1},z)=1_{M}$ (the identity operator on $M$).

(3) For ($\Z_{2}$-homogeneous) vectors $u,v\in V^{\bar 0}\cup V^{\bar 1}$, the following Jacobi
identity holds on $M$:
\begin{gather*}
z_{0}^{-1}\left(\frac{z_{1}-z_{2}}{z_{0}}\right)Y_{M}(u,z_{1})Y_{M}(v,z_{2})
-(-1)^{|u||v|}z_{0}^{-1}\delta\left(\frac{z_{2}-z_{1}}{-z_{0}}\right)Y_{M}(v,z_{2})Y_{M}(u,z_{1})\\
=z_{2}^{-1}\delta\left(\frac{z_{1}-z_{0}}{z_{2}}\right)Y_{M}(Y(u,z_{0})v,z_{2}).
\end{gather*}
} \end{defn}

Let $(M,Y_{M})$ be a weak $V$-module. Write $Y_{M}(\omega,z)=\sum_{n\in\mathbb{Z}}L(n)z^{-n-2}$,
recalling that $\omega$ is the conformal vector of $V$. The following
are consequences:
\begin{align}
 & [L(m),L(n)]=(m-n)L(m+n)+\delta_{m+n,0}\frac{m^{3}-m}{12}c_{V},\\
 & [L(-1),Y_{M}(v,z)]=Y_{M}(L(-1)v,z)=\frac{d}{dz}Y_{M}(v,z)\quad\text{for }v\in V,
\end{align}
where $c_{V}$ is the central charge of the Virasoro algebra acting
on $V$.

A {\em $\frac{1}{2}\mathbb{N}$-graded weak $V$-module} is a weak
$V$-module $M$ equipped with a $\frac{1}{2}\mathbb{N}$-grading
$M=\oplus_{n\in\frac{1}{2}\mathbb{N}}M(n)$ such that for $v\in V_{(h)}$
with $h\in\frac{1}{2}\Z$,
\begin{align}
v_{n}M(q)\subset M(h+q-n-1)\quad\text{ for }n\in\Z,\ q\in\frac{1}{2}\Z,
\end{align}
where $M(q)=0$ for $q\in-\frac{1}{2}\Z_{+}$ by convention.

A {\em $V$-module} is a weak $V$-module $M$ satisfying the condition
that $M=\oplus_{h\in\C}M_{(h)}$, where
\[
M_{(h)}=\{w\in M\ |\ L(0)w=hw\},
\]
$\dim M_{(h)}<\infty$ for all $h\in\C$, and for every fixed $h\in\C$,
$M_{(h+n)}=0$ for $n\in\frac{1}{2}\Z$ sufficiently negative.

Suppose that $M$ is a (nonzero) indecomposable $V$-module. Then there exists
a unique complex number $h$ such that $W=\oplus_{n\in\frac{1}{2}\N}W_{(h+n)}$
with $W_{(h)}\ne0$, and $W=\oplus_{n\in\frac{1}{2}\N}W(n)$ with
$W(n)=W_{(n+h)}$ is a $\frac{1}{2}\mathbb{N}$-graded weak $V$-module.
It follows (cf. \cite{DLM}) that every $V$-module can be made into
a $\frac{1}{2}\mathbb{N}$-graded weak $V$-module.

A {\em $\Z_{2}$-graded weak $V$-module} is a weak $V$-module
$W$ equipped with a $\Z_{2}$-grading $W=W^{\bar{0}}\oplus W^{\bar{1}}$
such that
\begin{align}
v_{n}W^{\beta}\subset W^{\alpha+\beta}\quad\text{ for }v\in V^{\alpha},\ \alpha,\beta\in\Z_{2},\ n\in\Z.
\end{align}
Let $W=\oplus_{q\in\frac{1}{2}\N}W(q)$ be a $\frac{1}{2}\N$-graded
weak $V$-module with $W(0)\ne0$. Set
\begin{align}
W^{\bar{0}}=\oplus_{n\in\N}W(n),\quad W^{\bar{1}}=\oplus_{n\in\N}W(n+1/2).
\end{align}
Then $W=W^{\bar{0}}\oplus W^{\bar{1}}$ is a $\Z_{2}$-graded weak
$V$-module. We shall always view a nonzero $\frac{1}{2}\N$-graded weak $V$-module
$W$ with $W(0)\ne0$ as a $\Z_{2}$-graded weak $V$-module in this
particular way.

\section{A family of vertex operator superalgebras}

In this section, we study a particular family of vertex operator superalgebras
which can be generated by the homogeneous subspaces of weight $2$ and weight  $\frac{3}{2}$.

We begin by recalling the following notion due to Kac (see \cite{Kac}):

\begin{defn}
{\em A confomal vertex superalgebra $V$ is said to
be {\em strongly generated} by a subset $T$ if $V$ is
linearly spanned by ${\bf 1}$ and the vectors of the form $a_{-n_{1}}^{(1)}\cdot\cdot\cdot a_{-n_{r}}^{(r)}{\bf 1}$
with $r\ge1$, $a^{(i)}\in T$ and $n_{i}\in\mathbb{Z}_{+}$.}
\end{defn}

The following is a preliminary result:

\begin{lemma}\label{basic1}
Let $V=\oplus_{q\in\frac{1}{2}\Z}V_{(q)}$ be a conformal vertex superalgebra
such that $V$ is  strongly generated by $V_{(\frac{3}{2})}+V_{(2)}$.
Then $V_{(q)}=0$ for $q\in -\frac{1}{2}\mathbb{Z}_{+}$,  $V_{(0)}=\C{\bf 1}$, $V_{(\frac{1}{2})}=0$,
$V_{(1)}=0$, and
\begin{eqnarray}
\quad V_{(\frac{5}{2})}=L(-1)V_{(\frac{3}{2})}.
\end{eqnarray}
\end{lemma}

\begin{proof} From assumption, $V$ is linearly spanned by ${\bf 1}$
and by the vectors $a_{-n_{1}}^{(1)}\cdot\cdot\cdot a_{-n_{r}}^{(r)}{\bf 1}$
for $r\ge1$, $a^{(i)}\in V_{(\frac{3}{2})}\cup V_{(2)},\  n_{i}\in\Z_{+}$.
As ${\rm wt}\;a^{(i)}=\frac{3}{2}$ or $2$, we have
\[
\text{wt}\;a_{-n_{i}}^{(i)}={\rm wt}\;a^{(i)}+n_{i}-1\ge\frac{3}{2}.
\]
Then $\wt(a_{-n_{1}}^{(1)}\cdots a_{-n_{r}}^{(r)}{\bf 1})\ge\frac{3}{2}r$.
It follows that $V_{(q)}=0$ for $q<0$, $V_{(0)}=\C{\bf 1}$,
$V_{(\frac{1}{2})}=0$, and $V_{(1)}=0$. It also follows that $V_{(\frac{5}{2})}$
is linearly spanned by $u_{-2}{\bf 1}$ for $u\in V_{(\frac{3}{2})}$.
As $v_{-2}{\bf 1}=L(-1)v$ for $v\in V$, we have $V_{(\frac{5}{2})}=L(-1)V_{(\frac{3}{2})}$.
\end{proof}

From now on, we assume that $V=\oplus_{q\in\frac{1}{2}\Z}V_{(q)}$ is a conformal vertex superalgebra
such that  $V_{(\frac{3}{2})}+V_{(2)}$ strongly generates $V$ unless it is stated otherwise.
Then
\begin{align}
V=V_{(0)}\oplus V_{(\frac{3}{2})}\oplus V_{(2)}\oplus V_{(\frac{5}{2})}\oplus V_{(3)}\oplus\cdots
\end{align}
with $V_{(0)}=\C {\bf 1}$.
As $L(1)V_{(2)}\subset V_{(1)}=0$ and $L(1)V_{(\frac{3}{2})}\subset V_{(\frac{1}{2})}=0$,
by (\ref{quasi-primary}) we immediately have:

\begin{coro}\label{L(1)a_0 and a_1, u_0 and u_1} The following relations
hold on $V$ for $a\in V_{(2)},\ u\in V_{(\frac{3}{2})},\ n\in \mathbb{Z}$:
\begin{eqnarray}\label{L(1)-brackets}
[L(1),a_{n}]=(2-n)a_{n+1},\ \ \ \ [L(1),u_{n}]=(1-n)u_{n+1}.
\end{eqnarray}
\end{coro}

The following are more technical results:

\begin{lemma}\label{a_0 v and a_1 v}
For $a\in V_{(2)},\ v\in V_{(\frac{3}{2})},$ we have
\begin{eqnarray}
a_{0}v=\frac{2}{3}L(-1)a_{1}v,\ \ \ \ v_{0}a=\frac{1}{3}L(-1)v_{1}a.
\end{eqnarray}
\end{lemma}

\begin{proof} Let $a\in V_{(2)}$, $v\in V_{(\frac{3}{2})}$. Since
$a_{0}v\in V_{(\frac{5}{2})}=L(-1)V_{(\frac{3}{2})}$ (by Lemma \ref{basic1}),
we have $a_{0}v=L(-1)u$ for some $u\in V_{(\frac{3}{2})}.$
On the one hand, we have
\begin{align*}
L(1)(L(-1)u)=L(-1)L(1)u+2L(0)u=3u,
\end{align*}
and  on the other hand, using (\ref{L(1)-brackets}) we have
\begin{align*}
L(1)(a_{0}v)=a_{0}L(1)v+[L(1),a_{0}]v=a_{0}L(1)v+2a_{1}v=2a_{1}v,
\end{align*}
noticing that $L(1)v\in V_{(\frac{1}{2})}=0$. Consequently, we get
$3u=2a_{1}v$. Thus $a_{0}v=\frac{2}{3}L(-1)a_{1}v$.
Similarly, writing $v_0a=L(-1)u$ for some $u\in V_{(\frac{3}{2})}$, we have
$L(1)L(-1)u=3a$ and $L(1)v_0a=v_1a$, which imply $v_{0}a=\frac{1}{3}L(-1)v_{1}a$.
\end{proof}

\begin{defn}\label{subspace-P}
{\em  Define
\begin{align}
P={\rm span}\left\{ u_{-1}v-\frac{1}{2}L(-1)u_0v\ \mid \ u,v\in V_{(\frac{3}{2})}\right\},
\end{align}
a subspace of $V_{(3)}$.}
\end{defn}

\begin{lemma}\label{P-direct-sum}
Let $V=\oplus_{q\in\frac{1}{2}\Z}V_{(q)}$
be a conformal vertex superalgebra such that $V_{(\frac{3}{2})}+V_{(2)}$ strongly generates $V$. Then
$L(1)P=0$ and
\begin{align}
V_{(3)}=L(-1)V_{(2)}\oplus P.
\end{align}
Furthermore, for any $a,b\in V_{(2)}$, there exists (uniquely) $f\in P$ such that
\begin{align}
a_0b=\frac{1}{2}L(-1)a_1b+f.
\end{align}
\end{lemma}

\begin{proof} For $u,v\in V_{(\frac{3}{2})}$, noticing that $L(1)v=0$ and $L(1)(u_0v)\in L(1)V_{(2)}=0$ we have
$$L(1)\!\left(u_{-1}v-\frac{1}{2}L(-1)u_0v\right)\!=u_{-1}L(1)v+2u_0v-\frac{1}{2}L(-1)L(1)u_0v-L(0)u_0v=2u_0v-2u_0v=0. $$
Thus $L(1)P=0$.
Since $V$ is strongly generated by $V_{(\frac{3}{2})}+V_{(2)}$, $V_{(3)}$ is linearly spanned by vectors
$a_{-2}{\bf 1}$ for $a\in V_{(2)}$ and $u_{-1}v$ $(=u_{-1}v_{-1}{\bf 1})$ for $u,v\in V_{(\frac{3}{2})}$.
Note that $a_{-2}{\bf 1}=L(-1)a$ for $a\in V_{(2)}$ and that $u_0v\in V_{(2)}$ for $u,v\in V_{(\frac{3}{2})}$.
Then we have $V_{(3)}=L(-1)V_{(2)}+P$.
This sum is a direct sum because for any $a\in V_{(2)}$, $L(-1)L(1)(L(-1)a)=4(L(-1)a)$, whereas $L(-1)L(1)P=0$.
Let $a,b\in V_{(2)}$. Then  $a_0b=L(-1)c+f$ for some $c\in V_{(2)},\ f\in P$. On the one hand, we have
$$L(1)(a_0b)=a_0L(1)b+2a_1b=2a_1b,$$
 and on the other hand, we have
$$L(1)(L(-1)c+f)=L(-1)L(1)c+2L(0)c+L(1)f=0+4c+0=4c.$$
Consequently, we get $c=\frac{1}{2}a_1b$ and then $a_0b=\frac{1}{2}L(-1)a_1b+f$.
\end{proof}

We next consider the conformal vertex subalgebra $\langle V_{(2)}\rangle$ generated by $V_{(2)}$.
Note that in view of Lemma \ref{P-direct-sum}, $V_{(3)}=L(-1)V_{(2)}$  if and only if
\begin{align}
 u_{-1}v=\frac{1}{2}L(-1)u_0v  \quad \text{  for }u,v\in V_{(\frac{3}{2})}.
 \end{align}

The following is essentially proved by Lam (see \cite{Lam}):

\begin{prop}\label{lam}
Let $V=\oplus_{n\in \mathbb{Z}}V_{(n)}$ be a conformal vertex algebra such that $V_{(2)}$ generates $V$.
Then all the following statements are equivalent:

(1) $V_{(2)}$ strongly generates $V$.

(2) $V_{(n)}=0$ for $n<0$, $V_{(0)}=\C {\bf 1}$, $V_{(1)}=0$, and  for $a,b\in V_{(2)}$, $a_0b=\frac{1}{2}L(-1)a_1b$.

(3)  For $a,b\in V_{(2)}$, $m,n\in \mathbb{Z}$,
\begin{align}\label{am-bn-bracket-1}
[a_m,b_n]=\frac{1}{2}(m-n)(a_{1}b)_{m+n-1}+\binom{m}{3}\langle a,b\rangle\delta_{m+n-2,0},
\end{align}
where $\langle a,b\rangle\in \C$ such that $a_3b=\langle a,b\rangle {\bf 1}$.
\end{prop}

\begin{proof}  Assume that $\langle V_{(2)}\rangle$ is strongly generated by $V_{(2)}$.
From the proof of Lemma \ref{P-direct-sum}, we see that $V_{(n)}=0$ for $n<0$,  $V_{(1)}=0$, $V_{(0)}=\C {\bf 1}$,
and that $a_0b\in L(-1)V_{(2)}$ for $a,b\in V_{(2)}$.
 Let $a,b\in V_{(2)}$. Then $a_{0}b=L(-1)c$ for some $c\in V_{(2)}$.
By the same reasoning as in the proof of Lemma \ref{P-direct-sum},  we have
$L(1)L(-1)c=4c$ and
$L(1)a_{0}b=a_{0}L(1)b+2a_{1}b=2a_{1}b$, which imply
$4c=2a_{1}b$, i.e.,  $c=\frac{1}{2}a_{1}b$.
Thus $a_{0}b=\frac{1}{2}L(-1)a_{1}b$. This shows that (1) implies (2).
Assume (2) holds true. Let $a,b\in V_{(2)},\ m,n\in \mathbb{Z}$.
As $a_3b\in V_{(0)}=\C {\bf 1}$, we have $a_3b=\langle a,b\rangle {\bf 1}$ with $\langle a,b\rangle\in \C$.
Then we get
\begin{align*}
[a_{m},b_{n}] & =(a_{0}b)_{m+n}+m(a_{1}b)_{m+n-1}+\binom{m}{2}(a_{2}b)_{m+n-2}+\binom{m}{3}(a_{3}b)_{m+n-3}\\
  & =\frac{1}{2}(L(-1)a_{1}b)_{m+n}+m(a_{1}b)_{m+n-1}+\binom{m}{3}\langle a,b\rangle \boldsymbol{1}_{m+n-3}\\
 & =-\frac{1}{2}(m+n)(a_{1}b)_{m+n-1}+m(a_{1}b)_{m+n-1}+\binom{m}{3}\langle a,b\rangle \delta_{m+n-2,0}\\
 &=\frac{1}{2}(m-n)(a_{1}b)_{m+n-1}+\binom{m}{3}\langle a,b\rangle\delta_{m+n-2,0},
\end{align*}
proving (\ref{am-bn-bracket-1}). Note that (3) implies that operators $a_n$ for $a\in V_{(2)},\ n\in \mathbb{Z}$
and the identity operator on $V$ linearly span a Lie subalgebra $L$ of $\mathfrak{gl}(V)$.
As $V_{(2)}$ generates $V$, we have $V=U(L){\bf 1}$. Since $a_n{\bf 1}=0$ for $a\in V_{(2)},\ n\ge 0$,
it follows from the P-B-W theorem that $V$ is strongly generated by $V_{(2)}$.
This proves the equivalence for all the three statements.
\end{proof}

As an immediate consequence of Proposition \ref{lam} we have:

\begin{coro}\label{a_0b and a_1 b}
Let $V=\oplus_{q\in\frac{1}{2}\Z}V_{(q)}$
be a conformal vertex superalgebra such that $V_{(\frac{3}{2})}+V_{(2)}$ strongly generates $V$.
In addition, assume that the vertex subalgebra $\langle V_{(2)}\rangle$ is
 strongly generated by $V_{(2)}$. Then
\begin{eqnarray}
a_{0}b=\frac{1}{2}L(-1)a_{1}b\ \text{ and }\ a_{1}b_{1}=b_{1}a_{1} \  \text{on }V
\end{eqnarray}
for $a,b\in V_{(2)}$.
\end{coro}

The following are some basic results about the homogeneous subspace $V_{(2)}$:

\begin{prop}\label{comm-assoc-algebra}
Let $V=\oplus_{q\in\frac{1}{2}\Z}V_{(q)}$ be a conformal vertex superalgebra
such that $V_{(\frac{3}{2})}+V_{(2)}$ strongly generates $V$.
Equip the subspace $V_{(2)}$ with a non-associative algebra structure with
\begin{align}\label{multiplication-V2}
ab=a_{1}b\quad  \text{ for }a,b\in V_{(2)}
\end{align}
and define a bilinear form $\langle\cdot,\cdot\rangle$ on $V_{(2)}$ by
\begin{align}\label{form-V2}
a_{3}b=\langle a,b\rangle\boldsymbol{1}\quad  \text{ for }a,b\in V_{(2)}.
\end{align}
Then $V_{(2)}$ is commutative and $\frac{1}{2}\omega$ is an identity, and the bilinear form $\langle\cdot,\cdot\rangle$ is symmetric.
Furthermore, if we in addition assume that  $a_0b\in L(-1)V_{(2)}$ for $a,b\in V_{(2)}$,
then the non-associative algebra $V_{(2)}$
is actually associative and  the bilinear form is associative
in the sense that $\langle ab,c\rangle=\langle a,bc\rangle$ for $a,b,c\in V_{(2)}$.
\end{prop}

\begin{proof}  Let $a,b\in V_{(2)}$, which are even. By skew symmetry
(\ref{skew-symmetry-component}), we have
\begin{align}\label{V2-commutativity}
a_{1}b=b_{1}a+\sum_{i\ge1}\frac{1}{i!}(-1)^{i}L(-1)^{i}\left(b_{1+i}a\right)=b_{1}a,
\end{align}
noticing that $b_{2}a\in V_{(1)}=0,\ b_3a\in V_{(0)}=\C{\bf 1}$ and $b_{1+i}a\in V_{(2-i)}=0$
for $i\ge 3$. This proves that the defined operation on $V_{(2)}$
is commutative. For $a\in V_{(2)}$,
we have $\frac{1}{2}\omega_{1}a=\frac{1}{2}L(0)a=a$. Thus $\frac{1}{2}\omega$ is an identity.

If we also assume that $a_0b\in L(-1)V_{(2)}$ for $a,b\in V_{(2)}$, by Corollary \ref{a_0b and a_1 b},
we have $a_{1}(b_{1}c)=b_{1}(a_{1}c)$. Then it follows from a classical
fact that $V_{(2)}$ is a commutative associative algebra with identity
$\frac{1}{2}\omega$.

Note that as $a_{3}b\in V_{(0)}=\mathbb{C}\boldsymbol{1}$ for $a,b\in V_{(2)}$,
the bilinear form $\langle \cdot,\cdot\rangle $ on $V_{(2)}$
is well defined. By skew symmetry we have
\[
a_{3}b=b_{3}a+\sum_{i\ge1}\frac{1}{i!}(-1)^{i}L(-1)^{i}\left(b_{3+i}a\right)=b_{3}a
\]
since $b_{3+i}a\in V_{(-i)}=0$ for $i\ge1$. Then
$\langle a,b\rangle =\langle b,a\rangle $.
This proves that the bilinear form is symmetric.

Now, assume $a_0b\in L(-1)V_{(2)}$ for all $a,b\in V_{(2)}$.
For $a,b,c\in V_{(2)}$, noticing that $b_{4}c\in V_{(-1)}=0$ and
$b_{3}c\in V_{(0)}=\C{\bf 1}$, using (\ref{asso commutator}) and
Lemma \ref{a_0b and a_1 b}, we get
\begin{align*}
 & \left(a_{1}b\right)_{3}c\\
 =\ & a_{1}\left(b_{3}c\right)+b_{4}\left(a_{0}c\right)-a_{0}\left(b_{4}c\right)-b_{3}\left(a_{1}c\right)\\
=\ & b_{4}\left(a_{0}c\right)-b_{3}\left(a_{1}c\right)\\
=\ & \frac{1}{2}b_{4}L(-1)\left(a_{1}c\right)-b_{3}\left(a_{1}c\right)\\
=\ & \frac{1}{2}L(-1)b_{4}\left(a_{1}c\right)-\frac{1}{2}\left[L(-1),b_{4}\right]\left(a_{1}c\right)-b_{3}\left(a_{1}c\right)\\
=\ & 2b_{3}(a_{1}c)-b_{3}(a_{1}c)\\
 =\ & b_{3}(a_{1}c).
\end{align*}
That is, $\langle ab,c\rangle =\langle b,ac\rangle $.
Consequently, we have $\langle ab,c\rangle=\langle ba,c\rangle=\langle a,bc\rangle$.
\end{proof}

Next, we consider the subspace $V_{(\frac{3}{2})}$. Note that for $u,v\in V_{(\frac{3}{2})}$,
as $\text{wt}\,u_{n}v=2-n$ for $n\in\Z$, we have
\begin{align}
u_{0}v\in V_{(2)},\ \ u_{1}v\in V_{(1)}=0,\ \ u_{2}v\in V_{(0)}=\mathbb{C}{\bf 1},\  \text{ and }u_nv=0\ \text{ for }n\ge 3.
\end{align}

\begin{defn}\label{forms-3/2}
{\em We equip $V_{(\frac{3}{2})}$
with a $\C$-valued bilinear form $\langle\cdot,\cdot\rangle$ and
a bilinear map
\begin{align}
\circ:\ \ V_{(\frac{3}{2})}\times V_{(\frac{3}{2})}\rightarrow V_{(2)};\quad(u,v)\mapsto u\circ v,
\end{align}
which are defined by
\begin{eqnarray}
u_{2}v=\langle u,v\rangle{\bf 1},\ \ \ \ u\circ v=u_{0}v\quad\text{ for }u,v\in V_{(\frac{3}{2})}.
\end{eqnarray}
} \end{defn}

Alternatively, we view the bilinear map $\circ$ as a $V_{(2)}$-valued
bilinear form on $V_{(\frac{3}{2})}$, denoted by $(\cdot,\cdot)$.

\begin{prop}\label{module}
Let $V=\oplus_{q\in\frac{1}{2}\Z}V_{(q)}$ be a conformal vertex superalgebra
such that $V$ is strongly generated by $V_{(\frac{3}{2})}+V_{(2)}$.
Define a left action of the non-associative algebra $V_{(2)}$ on $V_{(\frac{3}{2})}$ by
\begin{align}
av=\frac{4}{3}a_{1}v\quad\text{ for }a\in V_{(2)},\ v\in V_{(\frac{3}{2})}.
\end{align}
Then the identity $\frac{1}{2}\omega$ of $V_{(2)}$ acts as identity on $V_{(\frac{3}{2})}$ and
\begin{align}
a(bv)=(ab)v\quad \text{ for }a,b\in V_{(2)},\ v\in V_{(\frac{3}{2})}.
\end{align}
Furthermore, the bilinear form $\langle\cdot,\cdot\rangle$ on $V_{(\frac{3}{2})}$
and the bilinear map $\circ$ on $V_{(\frac{3}{2})}\times V_{(\frac{3}{2})}$ are symmetric,
and the following conditions hold for $a\in V_{(2)},\ u,v\in V_{(\frac{3}{2})}$:
\begin{eqnarray}\label{main-conditions}
\langle au,v\rangle=\langle u,av\rangle=\frac{4}{3}\langle a,u\circ v\rangle,\ \ \ \ (au,v)=a(u,v)=(u,av).
\end{eqnarray}
\end{prop}

\begin{proof}  For $v\in V_{(\frac{3}{2})}$, by definition we have
\begin{align*}
\frac{1}{2}\omega\cdot v=\frac{1}{2}\cdot\frac{4}{3}\omega_{1}v=\frac{2}{3}L(0)v=v.
\end{align*}
Let $a,b\in V_{(2)},$ $v\in V_{(\frac{3}{2})}$. Using the iterate
formula, as $b_{2}v\in V_{(\frac{1}{2})}=0$ we get
\begin{eqnarray*}
(a_{1}b)_{1}v=\sum_{i=0}^{\infty}(-1)^{i}\binom{1}{i}\left(a_{1-i}b_{1+i}+b_{2-i}a_{i}\right)v
=a_{1}b_{1}v+b_{2}a_{0}v-b_{1}a_{1}v.
\end{eqnarray*}
As $a_0v=\frac{2}{3}L(-1)a_{1}v$ from Lemma \ref{a_0 v and a_1 v}, we have
\[
b_{2}a_{0}v=\frac{2}{3}b_{2}L(-1)a_{1}v=\frac{2}{3}L(-1)b_{2}a_{1}v-\frac{2}{3}[L(-1),b_{2}]a_{1}v=\frac{4}{3}b_{1}a_{1}v,
\]
noticing that $b_{2}a_{1}v\in V_{(\frac{1}{2})}=0$. Combining these
two relations, we obtain
\begin{eqnarray}
(a_{1}b)_{1}v=a_{1}(b_{1}v)+\frac{4}{3}b_{1}(a_{1}v)-b_{1}(a_{1}v)=a_{1}(b_{1}v)+\frac{1}{3}b_{1}(a_{1}v).
\end{eqnarray}
Using symmetry we also have  $(b_{1}a)_{1}v=b_{1}(a_{1}v)+\frac{1}{3}a_{1}(b_{1}v)$.
As $a_{1}b=b_{1}a$, we obtain $a_{1}(b_{1}v)=b_{1}(a_{1}v)$. Consequently, we have
$(a_{1}b)_{1}v=\frac{4}{3}a_{1}(b_{1}v)$. This proves that $(ab)v=a(bv)$.

About the bilinear forms on $V_{(\frac{3}{2})}$, let $u,v\in V_{(\frac{3}{2})}$.
As $u,v$ are both odd, we have
\[
u_{n}v=(-1)^{n}\left(v_{n}u+\sum_{i\ge1}\frac{1}{i!}(-1)^{i}L(-1)^{i}(v_{n+i}u\right)
\]
for $n\in\mathbb{Z}$. In particular, we have
\[
u_{2}v=v_{2}u+\sum_{i\ge1}\frac{1}{i!}(-1)^{i+2}L(-1)^{i}\left(v_{2+i}u\right)=v_{2}u
\]
since $v_{2+i}u\in V_{(-i)}=0$ for $i\ge1$. Similarly, we have
\[
u_{0}v=v_{0}u-L(-1)v_{1}u+\frac{1}{2}L(-1)^{2}v_{2}u+\cdots=v_{0}u
\]
as $v_{1}u\in V_{(1)}=0$, $v_{2}u\in V_{(0)}=\C{\bf 1}$, and $v_{j}u\in V_{(2-j)}=0$
for $j\ge3$. Thus both the $\C$-valued bilinear form $\langle\cdot,\cdot\rangle$
and the $V_{(2)}$-valued bilinear form $(\cdot,\cdot)$ on $V_{(\frac{3}{2})}$
are symmetric.

Furthermore, for $a\in V_{(2)},\ u,v\in V_{(\frac{3}{2})}$, we have
\begin{eqnarray}
(a_{1}u)_{2}v=a_{1}u_{2}v+u_{3}a_{0}v-a_{0}u_{3}v-u_{2}a_{1}v=u_{3}a_{0}v-u_{2}a_{1}v,
\end{eqnarray}
noticing that $u_{3}v\in V_{(-1)}=0$ and $a_{1}u_{2}v=0$ as $u_{2}v\in V_{(0)}=\C{\bf 1}$.
As $a_{0}v=\frac{2}{3}L(-1)a_{1}v$ by Lemma \ref{a_0 v and a_1 v},
we have
\[
u_{3}a_{0}v=\frac{2}{3}u_{3}L(-1)a_{1}v=\frac{2}{3}L(-1)u_{3}a_{1}v+2u_{2}a_{1}v=2u_{2}a_{1}v.
\]
Thus $(a_{1}u)_{2}v=u_{2}a_{1}v$. This proves $\langle au,v\rangle=\langle u,av\rangle$.
Noticing that $a_{3}u=0=a_{3}v$, $a_{2}v=0$, and $a_{0}u=\frac{2}{3}L(-1)a_{1}u$
we obtain
\begin{eqnarray*}
 &  & \langle a,u\circ v\rangle=a_{3}u_{0}v=[a_{3},u_{0}]v=(a_{0}u)_{3}v+3(a_{1}u)_{2}v+3(a_{2}u)_{1}v+(a_{3}u)_{0}v\\
 & = & \frac{2}{3}(L(-1)a_{1}u)_{3}v+3(a_{1}u)_{2}v=-2(a_{1}u)_{2}v+3(a_{1}u)_{2}v=(a_{1}u)_{2}v=\frac{3}{4}\langle au,v\rangle,
\end{eqnarray*}
as desired.

On the other hand, we have
\begin{eqnarray}
(a_{1}u)_{0}v=a_{1}u_{0}v+u_{1}a_{0}v-a_{0}u_{1}v-u_{0}a_{1}v=a_{1}u_{0}v+u_{1}a_{0}v-u_{0}a_{1}v
\end{eqnarray}
as $u_{1}v\in V_{(1)}=0$. With $a_{0}v=\frac{2}{3}L(-1)a_{1}v$ again,
we have
\[
u_{1}a_{0}v=\frac{2}{3}u_{1}L(-1)a_{1}v=\frac{2}{3}L(-1)(u_{1}a_{1}v)+\frac{2}{3}u_{0}a_{1}v=\frac{2}{3}u_{0}a_{1}v
\]
because $u_{1}a_{1}v\in V_{(1)}=0$. Then we obtain
\begin{eqnarray}
(a_{1}u)_{0}v=a_{1}u_{0}v-\frac{1}{3}u_{0}a_{1}v.
\end{eqnarray}
That is, $(au,v)=\frac{4}{3}a(u,v)-\frac{1}{3}(u,av)$ (as $aw=\frac{4}{3}a_{1}w$
for $w\in V_{(\frac{3}{2})}$). Switching $u$ and $v$ we get $(av,u)=\frac{4}{3}a(v,u)-\frac{1}{3}(v,au)$.
As $(\cdot,\cdot)$ is symmetric, we have $(u,av)=\frac{4}{3}a(u,v)-\frac{1}{3}(au,v)$.
Then by solving the two equations we obtain $(au,v)=(u,av)=a(u,v)$.
\end{proof}

Furthermore, we have:

\begin{prop}\label{vosa-commutator-relations}
Let $V=\oplus_{q\in\frac{1}{2}\Z}V_{(q)}$ be a conformal vertex superalgebra
such that  $V_{(\frac{3}{2})}+V_{(2)}$ strongly generates $V$ and such that $\langle V_{(2)}\rangle$ is
strongly generated by $V_{(2)}$.  Then
\begin{eqnarray}
 &  & [a_{m},b_{n}]=\frac{1}{2}(m-n)(a_{1}b)_{m+n-1}+\binom{m}{3}\langle a,b\rangle\delta_{m+n-2,0},\label{am-bn-last}\ \ \ \ \\
 &  & [a_{m},u_{n}]=\frac{1}{3}(m-2n)(a_{1}u)_{m+n-1},\label{am-un-bracket}\\
 &  & [u_{m},v_{n}]_{+}=(u_{0}v)_{m+n}+\frac{1}{2}m(m-1)\langle u,v\rangle\delta_{m+n-1,0}
\end{eqnarray}
for $a,b\in V_{(2)}$, $u,v\in V_{(\frac{3}{2})}$, $m,n\in\mathbb{Z}$, where $[u_{m},v_{n}]_{+}=u_{m}v_{n}+v_{n}u_{m}$.
\end{prop}

\begin{proof} The first relation has been obtained in Proposition \ref{lam} (see (\ref{am-bn-bracket-1})).
 For the other two bracket relations, noticing that
$a_{i}u\in V_{(\frac{5}{2}-i)}=0$ for $i\ge2$ and $u_{i}v\in V_{(2-i)}=0$
for $i\ge3$, then using Lemma \ref{a_0 v and a_1 v} we get
\begin{align*}
[a_{m},u_{n}] & =(a_{0}u)_{m+n}+m(a_{1}u)_{m+n-1}\\
 & =\frac{2}{3}(L(-1)a_{1}u)_{m+n}+m(a_{1}u)_{m+n-1}\\
 & =-\frac{2}{3}(m+n)(a_{1}u)_{m+n-1}+m(a_{1}u)_{m+n-1}\\
 & =\frac{1}{3}(m-2n)(a_{1}u)_{m+n-1},
\end{align*}
and
\begin{align*}
[u_{m},v_{n}]_{+}=(u_{0}v)_{m+n}+\binom{m}{2}(u_{2}v)_{m+n-2}
=(u_{0}v)_{m+n}+\frac{1}{2}m(m-1)\langle u,v\rangle\delta_{m+n-1,0},
\end{align*}
where we also use the fact that $u_{1}v\in V_{(1)}=0$.
\end{proof}

\begin{rem}\label{bracket-equivalence}
{\em Let $V=\oplus_{q\in\frac{1}{2}\Z}V_{(q)}$ be a general conformal vertex superalgebra
such that  $V_{(\frac{3}{2})}+V_{(2)}$ strongly generates $V$. Here, we show that for $a,b\in V_{(2)}$,
(\ref{am-bn-bracket-1}) (or (\ref{am-bn-last})) holds for all $m,n\in \Z$ if and only if $a_0b=\frac{1}{2}L(-1)a_1b$.
The ``if'' part follows from the proof of Proposition \ref{lam}.
Assume (\ref{am-bn-bracket-1}) holds. On the other hand, Borcherds' commutator formula gives
\begin{align*}
[a_m,b_n]=(a_0b)_{m+n}+m(a_1b)_{m+n-1}+\binom{m}{3}\langle a,b\rangle\delta_{m+n-2,0}.
\end{align*}
Consequently, we have $(a_0b)_{m+n}=-\frac{1}{2}(m+n)(a_1b)_{m+n-1}$.
Applying this to ${\bf 1}$ and taking $m=0,\ n=-1$, we get
$(a_0b)_{-1}{\bf 1}=\frac{1}{2}(a_1b)_{-2}{\bf 1}$. That is, $a_0b=\frac{1}{2}L(-1)a_1b$.}
\end{rem}

\begin{rem}\label{voa-Lie-algebra-affinization}
{\em Let $V$ be a general vertex superalgebra (with only a $\Z_{2}$-grading).
Set $L(V)=V\otimes\C[t,t^{-1}]$ and $\widehat{L}(-1)=L(-1)\otimes1+1\otimes\frac{d}{dt}$.
Then (see \cite{Bor}, \cite{FFR}, \cite{Li-thesis}) we have a Lie
superalgebra structure on $L(V)/\widehat{L}(-1)L(V)$ with
\begin{eqnarray}
[u(m),v(n)]_{\pm}=\sum_{i\ge0}\binom{m}{i}(u_{i}v)(m+n-i)
\end{eqnarray}
for $\Z_{2}$-homogeneous vectors $u,v\in V$ and for $m,n\in\Z$,
where $a(k)$ denotes the image of $a\otimes t^{k}$ in $L(V)/\widehat{L}(-1)L(V)$
for $a\in V,\ k\in\Z$.}
\end{rem}

\begin{coro}\label{Lie-subalgebra}
Let $V=\oplus_{q\in\frac{1}{2}\Z}V_{(q)}$
be a conformal vertex superalgebra such that $V$ is strongly generated by
the subspace $V_{(\frac{3}{2})}+V_{(2)}$ and such that $a_0b\in L(-1)V_{(2)}$ for $a,b\in V_{(2)}$.
Then the vector space $(V_{(2)}+V_{(\frac{3}{2})})\otimes\C[t,t^{-1}]+\C K$
is a Lie superalgebra with the Lie bracket relations given as in Proposition
\ref{vosa-commutator-relations}.
\end{coro}

\begin{proof} From Remark \ref{voa-Lie-algebra-affinization}, we
have a Lie superalgebra $L(V)/\widehat{L}(-1)L(V)$. Define a linear map
\[
\psi:(V_{(2)}+V_{(\frac{3}{2})})\otimes\C[t,t^{-1}]+\C K\rightarrow L(V)/\widehat{L}(-1)L(V)
\]
by
\begin{align*}
\psi(K)={\bf 1}(-1),\quad\psi(v\otimes t^{n})=v(n)\ \ \text{ for }v\in V_{(2)}+V_{(\frac{3}{2})},\ n\in\Z.
\end{align*}
It is clear that the image of $\psi$ in $L(V)/\widehat{L}(-1)L(V)$
is a Lie sub-superalgebra. Notice that
\[
\widehat{L}(-1)L(V)\cap\left((V_{(2)}+V_{(\frac{3}{2})})\otimes\C[t,t^{-1}]+\C K\right)=0.
\]
Then $\psi$ is one-to-one. It follows that $(V_{(2)}+V_{(\frac{3}{2})})\otimes\C[t,t^{-1}]+\C K$
is a Lie superalgebra with the desired bracket relations.
\end{proof}

Next, we consider the general case where $\langle V_{(2)}\rangle$ is  {\em not necessarily} strongly generated by $V_{(2)}$.
In the general case, $A=V_{(2)}$ is not necessarily associative.
Recall that  $V_{(3)}=L(-1)V_{(2)}$ if and only if $P=0$, where
$$P=\text{span}\left\{ u_{-1}v-\frac{1}{2}L(-1)u_0v\ \mid \ u,v\in V_{(\frac{3}{2})}\right\}.$$
Also, recall that for $u,v\in V_{(\frac{3}{2})}$, we have
\begin{align*}
  u_1v=0,\  \   u_2v=\langle u,v\rangle {\bf 1},\   \  u_rv=0\  \text{ for }r\ge 3,\
\end{align*}
which (by skew symmetry) imply
\begin{align*}
u_0v=v_0u,\  \  u_{-1}v=-v_{-1}u+L(-1)v_0u.
\end{align*}
Then we get
\begin{align}\label{skew-uv}
u_{-1}v-\frac{1}{2}L(-1)u_0v=-\left(v_{-1}u-\frac{1}{2}L(-1)v_0u\right)\!.
\end{align}

\begin{rem}\label{dim-U=1}
{\em Assume $\dim V_{(\frac{3}{2})}=1$. From (\ref{skew-uv}) we have
$u_{-1}v-\frac{1}{2}L(-1)u_0v=0$ for $u,v\in V_{(\frac{3}{2})}$, i.e.,  $P=0$. Thus $V_{(3)}=L(-1)V_{(2)}$.
In this case,  $V_{(2)}$ is a commutative associative algebra by Proposition \ref{comm-assoc-algebra}.}
\end{rem}

\begin{rem}\label{dim-U>1}
{\em Assume that the bilinear form on $V_{(\frac{3}{2})}$ is non-degenerate.
For $u,v,w\in V_{(\frac{3}{2})}$, we have
\begin{align}\label{u2-v-w-calculation}
&u_2\left(v_{-1}w-\frac{1}{2}L(-1)v_0w\right)\nonumber\\
= &-v_{-1}u_2w+ (u_0v)_{1}w+(u_2v)_{-1}w-\frac{1}{2}L(-1)u_2(v_0w)-u_{1}v_0w\nonumber\\
= &-v_{-1}u_2w+ (u_0v)_{1}w+(u_2v)_{-1}w-0+v_0u_1w-(u_0v)_{1}w\nonumber\\
= &-v_{-1}u_2w+(u_2v)_{-1}w\nonumber\\
= & -\langle u,w\rangle v+\langle u,v\rangle w,
\end{align}
which is not zero if $v$ and $w$ are linearly independent and if either $\langle u,w\rangle\ne 0$ or $\langle u,v\rangle\ne 0$.
From this, we see that  if $\dim V_{(\frac{3}{2})}>1$, then $P\ne 0$, or equivalently, $V_{(3)}\ne L(-1)V_{(2)}$.}
\end{rem}

\begin{rem}\label{dim-U=2}
{\em Consider the case in which $\dim V_{(\frac{3}{2})}=2$ and the bilinear form on $V_{(\frac{3}{2})}$ is non-degenerate.
Let $\{ u,v\}$ be an orthonormal basis of $V_{(\frac{3}{2})}$. We see that $P$
is $1$-dimensional and spanned by $u_{-1}v-\frac{1}{2}L(-1)u_0v$.
Then there is a skew symmetric bilinear form  $\mu(\cdot,\cdot)$ on $V_{(2)}$ such that
\begin{align}
a_0b=\frac{1}{2}L(-1)a_1b+\mu(a,b)\left(u_{-1}v-\frac{1}{2}L(-1)u_0v\right)\!
\end{align}
for $a,b\in V_{(2)}$. The vertex operator superalgebra $V$ is uniquely determined by a commutative non-associative algebra $V_{(2)}$,
a module $V_{(\frac{3}{2})}$, and a skew symmetric bilinear form $\mu(\cdot,\cdot)$ on $V_{(2)}$.}
\end{rem}

Let $a,b\in V_{(2)}$. Note that $a_2b=0$, $a_3b\in \C {\bf 1}$, $a_1b=b_1a$, and $a_0b=-b_0a+L(-1)b_1a$. Then
we have
\begin{align}
a_0b-\frac{1}{2}L(-1)a_1b=-\left(b_0a-\frac{1}{2}L(-1)b_1a\right).
\end{align}
In particular, we have $a_0a=\frac{1}{2}L(-1)a_1a$.

\begin{defn}
{\em Let $V$ be a vertex operator superalgebra.
A bilinear form $\langle\cdot,\cdot\rangle$ on $V$ is said to be {\em invariant} (see \cite{FHL}) if
\begin{align}
\langle Y(u,z)v,w\rangle=\left\langle v,Y\left(e^{zL(1)}z^{-2L(0)}e^{\pi i L(0)}v,z^{-1}\right)w\right\rangle
\end{align}
for $u,v,w\in V$, where $i=\sqrt{-1}$ and by definition $e^{\pi i L(0)}v=e^{q\pi i}v$
for $v\in V_{(q)}$ with $q\in \frac{1}{2}\mathbb{Z}$.}
\end{defn}

From  \cite{FHL},  any invariant bilinear form on a vertex operator superalgebra is automatically symmetric.
Suppose that $\langle\cdot,\cdot\rangle$  is an invariant bilinear form on
a vertex operator superalgebra $V=\oplus_{q\in\frac{1}{2}\Z}V_{(q)}$. Then (see  \cite{FHL})
\begin{align}
\langle L(n)u,v\rangle =\langle u,L(-n)v\rangle \quad \text{ for }n\in \mathbb{Z},\ u,v\in V.
\end{align}
In particular, $\langle L(0)u,v\rangle =\langle u,L(0)v\rangle$.
It follows that $\langle V_{(p)},V_{(q)}\rangle =0$ for $p,q\in \frac{1}{2}\mathbb{Z}$ with $p\ne q$.

\begin{prop}\label{Inv-form-V}
Let $V=\oplus_{q\in\frac{1}{2}\Z}V_{(q)}$
be a vertex operator superalgebra such that $V_{(\frac{3}{2})}+V_{(2)}$ strongly generates $V$.
Then $V$ admits a symmetric invariant bilinear form $\langle\cdot,\cdot\rangle$
which is uniquely determined by $\langle {\bf 1},{\bf 1}\rangle=1$.
Furthermore, the kernel of $\langle\cdot,\cdot\rangle$ coincides with the unique maximal ideal $I$ (with ${\bf 1}\notin I$).
On the other hand, $V$ is simple if and only if the bilinear form $\langle\cdot,\cdot\rangle$ is non-degenerate.
\end{prop}

\begin{proof} Recall that $V_{(1)}=0$ and $V_{(0)}=\C {\bf 1}$, in particular, $L(1)V_{(1)}=0$.
It then follows from \cite{Li-form, S}  that $V$ admits a symmetric invariant bilinear form $\langle\cdot,\cdot\rangle$
which is uniquely determined by $\langle {\bf 1},{\bf 1}\rangle=1$.
Note that the kernel of $\langle\cdot,\cdot\rangle$ is a proper ideal of $V$. If $V$ is simple,
then it follows that $\langle\cdot,\cdot\rangle$ is non-degenerate.

Let $I$ be any proper ideal of $V$. Note that $I$ is automatically graded.
Since $V_{(0)}=\C {\bf 1}$ and $\langle {\bf 1},{\bf 1}\rangle =1$, we conclude that $I_{(0)}=0$.
It follows that $\langle {\bf 1},v\rangle=0$ for all $v\in I$.
Then for any $v\in I,\ f\in V$, we have
$$\langle Y(f,z){\bf 1},v\rangle=\left\langle {\bf 1}, Y\left(e^{zL(1)}z^{-2L(0)}e^{\pi i L(0)}f,z^{-1}\right)v  \right\rangle=0,$$
which implies $\langle f,v\rangle=0$.  Thus $I$ is contained in the kernel of $\langle\cdot,\cdot\rangle$. It then follows that
the kernel of $\langle\cdot,\cdot\rangle$ is the unique maximal ideal $I$ with ${\bf 1}\notin I$. The last assertion is clear.
\end{proof}

Let $V=\oplus_{q\in\frac{1}{2}\Z}V_{(q)}$ be a vertex operator superalgebra such that
$V_{(\frac{3}{2})}+V_{(2)}$ strongly generates $V$, and let $\langle\cdot,\cdot\rangle_V$ be the
unique symmetric invariant bilinear form such that $\langle {\bf 1},{\bf 1}\rangle_V=1$.
For $a\in V_{(2)}$, as $L(1)a=0$ and $L(0)a=2a$, we have
\begin{align}
Y\left(e^{zL(1)}z^{-2L(0)}e^{\pi i L(0)}a,z^{-1}\right)=z^{-4}Y(a,z^{-1})=\sum_{m\in \mathbb{Z}}a_mz^{m-3}.
\end{align}
Then
\begin{align}\label{invariance-an}
\langle a_nu,v\rangle_V=\langle u, a_{2-n}v\rangle_V\quad \text{ for }u,v\in V,\ n\in \mathbb{Z}.
\end{align}
In particular, for $a,b\in V_{(2)}$ we have
\begin{align}
\langle a,b\rangle_V=\langle a_{-1}{\bf 1},b\rangle_V=\langle {\bf 1},a_{3}b\rangle_V
= \langle a,b\rangle\langle {\bf 1},{\bf 1}\rangle_V=\langle a,b\rangle,
\end{align}
where $\langle \cdot,\cdot\rangle$ is the bilinear form on $V_{(2)}$ defined before by
$a_3b=\langle a,b\rangle{\bf 1}$.
Similarly, for $u\in V_{(\frac{3}{2})}$ we have
\begin{align}
Y\left(e^{zL(1)}z^{-2L(0)}e^{\pi i L(0)}u,z^{-1}\right)=(-i)z^{-3}Y(u,z^{-1})=(-i)\sum_{m\in \mathbb{Z}}u_mz^{m-2}
\end{align}
and
\begin{align}
\langle u,v\rangle_V=\langle u_{-1}{\bf 1},v\rangle_V=(-i)\langle {\bf 1},u_{2}v\rangle_V
=(-i) \langle u,v\rangle\langle {\bf 1},{\bf 1}\rangle_V=(-i)\langle u,v\rangle.
\end{align}
As a consequence of Proposition \ref{Inv-form-V}, we immediately have:

\begin{coro}\label{bilinear-forms-A-U}
Suppose that $V=\oplus_{q\in\frac{1}{2}\Z}V_{(q)}$ is a simple vertex operator superalgebra such that
$V_{(\frac{3}{2})}+V_{(2)}$ strongly generates $V$. Then the bilinear form $\langle \cdot,\cdot\rangle$ on $V_{(2)}$ and
the bilinear form $\langle \cdot,\cdot\rangle$  on $V_{(\frac{3}{2})}$ are non-degenerate.
\end{coro}

We also have:

\begin{prop}\label{simple-key}
Suppose that $V=\oplus_{q\in\frac{1}{2}\Z}V_{(q)}$ is a simple vertex operator superalgebra such that
$V_{(\frac{3}{2})}+V_{(2)}$ strongly generates $V$. Then
\begin{align}
a_0b=\frac{1}{2}L(-1)a_1b\quad \text{ for }a,b\in V_{(2)}.
\end{align}
\end{prop}

\begin{proof} First, by Proposition \ref{Inv-form-V} there is
a (unique) non-degenerate invariant bilinear form $\langle\cdot,\cdot\rangle$ on $V$ with $\langle {\bf 1},{\bf 1}\rangle =1$.
In particular, $\langle \cdot,\cdot\rangle$ is non-degenerate on $V_{(3)}$.
Second, by Lemma \ref{P-direct-sum} we have $V_{(3)}=L(-1)V_{(2)}+P$ with $L(1)P=0$, where $P$ is the linear span
of vectors $u_{-1}v-\frac{1}{2}L(-1)u_0v$ for $u,v\in V_{(\frac{3}{2})}$. Now, let $a,b\in V_{(2)}$.
 By Lemma \ref{P-direct-sum} we have
$\Delta(a,b):=a_0b-\frac{1}{2}L(-1)a_1b\in P\subset V_{(3)}$. In particular, $L(1)\Delta(a,b)=0$.
Then
\begin{align}
\left\langle \Delta(a,b), L(-1)V_{(2)}\right\rangle=\left\langle L(1)\Delta(a,b),V_{(2)}\right\rangle=0.
\end{align}
Next, we prove
$\langle \Delta(a,b),P\rangle=0$, so that $\langle \Delta(a,b),V_{(3)}\rangle=0$. It will then follow that
$\Delta(a,b)=0$.

Let $u,v\in V_{(\frac{3}{2})}$. We have $L(1)\left(u_{-1}v-\frac{1}{2}L(-1)u_0v\right)=0$. On the other hand,
noticing that $a_2v\in V_{(\frac{1}{2})}=0$ and $(a_1u)_1v\in V_{(1)}=0$, using
Lemma \ref{a_0 v and a_1 v} and (\ref{main-conditions})  we have
\begin{align*}
a_2\left(u_{-1}v-\frac{1}{2}L(-1)u_0v\right)=\ &\frac{4}{3}(a_1u)_{0}v-\frac{1}{3}L(-1)(a_1u)_{1}v
-a_{1}u_0v\\
=\ &\frac{4}{3}(a_1u)_{0}v-a_{1}u_0v\\
=\ &\frac{4}{3}(a_1u)_{0}v-u_0a_{1}v-\frac{1}{3}(a_1u)_0v\\
=\ &(a_1u)_{0}v-u_0(a_{1}v)\\
=\ &0.
\end{align*}
Then using invariance property (\ref{invariance-an}) we obtain
\begin{align*}
&\left\langle \Delta(a,b), u_{-1}v-\frac{1}{2}L(-1)u_0v\right\rangle\\
=\ & \left\langle a_0b, u_{-1}v-\frac{1}{2}L(-1)u_0v\right\rangle-\left\langle \frac{1}{2}L(-1)a_1b, u_{-1}v-\frac{1}{2}L(-1)u_0v\right\rangle\\
=\ & \left\langle b, a_2\left(u_{-1}v-\frac{1}{2}L(-1)u_0v\right)\right\rangle
-\frac{1}{2}\left\langle a_1b, L(1)\left(u_{-1}v-\frac{1}{2}L(-1)u_0v\right)\right\rangle\\
=\ &\langle b,0\rangle-\frac{1}{2}\langle a_1b,0\rangle\\
=\ & 0.
\end{align*}
This proves $\langle \Delta(a,b),P\rangle=0$, completing the proof.
\end{proof}

\begin{rem}\label{L-kernel}
{\em Let $V=\oplus_{q\in\frac{1}{2}\Z}V_{(q)}$ be a vertex operator superalgebra such that
$V_{(\frac{3}{2})}+V_{(2)}$ strongly generates $V$. Set
\begin{align}
L=\text{span}\left\{ a_0b-\frac{1}{2}L(-1)a_1b\mid  a,b\in V_{(2)}\right\}\!,
\end{align}
a subspace of $V_{(3)}$. From the proof of Proposition \ref{simple-key}, we see that
$\langle L,V_{(3)}\rangle=0$. Then $L$ is contained in the kernel of the invariant bilinear form on $V$.}
\end{rem}

\section{Lie superalgebra $\mathcal{L}(A,U)$}

To establish the existence of the vertex operator superalgebras considered in Section 3,
in this section we construct a Lie superalgebra $\mathcal{L}(A,U)$
from a general commutative associative algebra $A$ and an $A$-module $U$, equipped with other structures,
satisfying the conditions (\ref{main-conditions}).  We shall construct the desired vertex operator superalgebras in the next section.

The following is the main result of this section:

\begin{thm}\label{iff-Lie-super-algebra}
Let $A$ be a non-associative algebra equipped with a bilinear form $\langle\cdot,\cdot\rangle$
and let $U$ be a vector space on which $A$ acts from left. Assume
that $U$ is equipped with a bilinear form $\langle\cdot,\cdot\rangle$
and with a bilinear map $\circ:U\times U\rightarrow A;\ (u,v)\mapsto u\circ v$.
Set
\begin{eqnarray}
\mathcal{L}(A,U)=A\otimes\C[t,t^{-1}]\oplus U\otimes\C[t,t^{-1}]\oplus\C K,
\end{eqnarray}
where $K$ is a nonzero element. View $\mathcal{L}(A,U)$ as a superspace
with
\begin{align}
\mathcal{L}(A,U)^{0}=A\otimes\C[t,t^{-1}]+\C K,\quad\mathcal{L}(A,U)^{1}=U\otimes\C[t,t^{-1}].
\end{align}
Define a (bilinear) operation $[\cdot,\cdot]$ on $\mathcal{L}(A,U)$
by
\begin{eqnarray}
 &  & \left[a\otimes t^{m},b\otimes t^{n}\right]=\frac{1}{2}(m-n)(ab)\otimes t^{m+n-1}+\binom{m}{3}\langle a,b\rangle\delta_{m+n-2,0}K,\label{am-bn-lie} \\
 &  & [a\otimes t^{m},u\otimes t^{n}]=\frac{1}{4}(m-2n)(au)\otimes t^{m+n-1}, \label{am-un-lie}\\
 &  & [u\otimes t^{n},a\otimes t^{m}]=-[a\otimes t^{m},u\otimes t^{n}],\nonumber \\
 &  & [u\otimes t^{m},v\otimes t^{n}]_{+}=(u\circ v)\otimes t^{m+n}+\frac{1}{2}m(m-1)\langle u,v\rangle\delta_{m+n-1,0}K,\nonumber \\
 &  & [a\otimes t^{m},K]=0=[u\otimes t^{m},K]
\end{eqnarray}
for $a,b\in A$, $u,v\in U$ and for $m,n\in\mathbb{Z}$. Then $(\mathcal{L}(A,U),[\cdot,\cdot])$
carries the structure of a Lie superalgebra if and only if the following
conditions are satisfied:

(i) $A$ is a commutative associative algebra and the bilinear form
$\langle\cdot,\cdot\rangle$ on $A$ is symmetric and associative.

(ii) $U$ is an $A$-module and both the bilinear form $\langle\cdot,\cdot\rangle$
on $U$ and the bilinear map $\circ$ are symmetric.

(iii) For $a\in A,\ u,v,w\in U$,
\begin{eqnarray}
 &  & \langle v,au\rangle=\langle u,av\rangle=\frac{4}{3}\langle a,u\circ v\rangle,\label{NS-algebra}\\
 &  & a(u\circ v)=v\circ au=u\circ av,\label{invariant-properties}\\
 &  & (u\circ v)w=(v\circ w)u=(w\circ u)v.
\end{eqnarray}
\end{thm}

\begin{proof} First of all, it can be readily seen that the Lie superalgebra
skew symmetry is equivalent to that $A$ is a commutative algebra
and the bilinear forms on $A$ and $U$ are symmetric. In view of
this, we now assume these properties for the rest of the proof.

Case 1: Consider $a,b,c\in A$. For $m,n,k\in\Z$, we have
\begin{align*}
 & \left[a\otimes t^{m},\left[b\otimes t^{n},c\otimes t^{k}\right]\right]\\
=\  & \left[a\otimes t^{m},\frac{1}{2}(n-k)(bc\otimes t^{n+k-1})+\binom{n}{3}\langle b,c\rangle\delta_{n+k-2,0}K\right]\\
=\  & \frac{1}{4}(n-k)(m-n-k+1)a(bc)\otimes t^{m+n+k-2}+\frac{1}{2}\binom{m}{3}(n-k)\langle a,bc\rangle\delta_{m+n+k-3,0}K.
\end{align*}
Using permutations we also have
\begin{align*}
 & \left[b\otimes t^{n},[c\otimes t^{k},a\otimes t^{m}]\right]\\
=\  & \frac{1}{4}(k-m)(n-k-m+1)b(ca)\otimes t^{m+n+k-2}+\frac{1}{2}\binom{n}{3}(k-m)\langle b,ca\rangle\delta_{m+n+k-3,0}K
\end{align*}
and
\begin{align*}
 & \left[c\otimes t^{k},[a\otimes t^{m},b\otimes t^{n}]\right]\\
=\  & \frac{1}{4}(m-n)(k-m-n+1)c(ab)\otimes t^{m+n+k-2}+\frac{1}{2}\binom{k}{3}(m-n)\langle c,ab\rangle\delta_{m+n+k-3,0}K.
\end{align*}
Then
\[
\left[a\otimes t^{m},[b\otimes t^{n},c\otimes t^{k}]\right]+\left[b\otimes t^{n},[c\otimes t^{k},a\otimes t^{m}]\right]+\left[c\otimes t^{k},[a\otimes t^{m},b\otimes t^{n}]\right]=0
\]
if and only if
\begin{align}
(n-k)(m-n-k+1)a(bc)+(k-m)(n-k-m+1)b(ca)+(m-n)(k-m-n+1)c(ab)=0\label{A-abc}
\end{align}
and
\begin{align}
\left(\binom{m}{3}(n-k)\langle a,bc\rangle+\binom{n}{3}(k-m)\langle b,ca\rangle+\binom{k}{3}(m-n)\langle c,ab\rangle\right)\delta_{m+n+k-3,0}=0.\label{A-abc-K}
\end{align}
Note that taking $m=n=k+1=1$ and $m=k=n+1=1$ in (\ref{A-abc}),
we get $a(bc)=b(ca)$ and $a(bc)=c(ab)$, respectively. With $A$
a commutative algebra, (\ref{A-abc}) implies that $A$ is commutative
and associative. Conversely, if $A$ is commutative and associative,
it is straightforward to show that (\ref{A-abc}) holds.

On the other hand, we see that (\ref{A-abc-K}) with $k=m=3$ and
$k=n=3$ implies $\langle a,bc\rangle=\langle c,ab\rangle$ and $\langle c,ab\rangle=\langle b,ca\rangle$,
respectively. 
Conversely, (\ref{A-abc-K}) also follows from $\langle a,bc\rangle=\langle c,ab\rangle$
and $\langle c,ab\rangle=\langle b,ca\rangle$ as
\begin{align}
\left((n-k)\binom{m}{3}+(k-m)\binom{n}{3}+(m-n)\binom{k}{3}\right)\delta_{m+n+k,3}=0,\label{binom-sum=00003D00003D0}
\end{align}
which can be proved straightforwardly.

Case 2: Consider $a,b\in A,\ u\in U$. Let $m,n,k\in\mathbb{Z}$.
By definition, we have
\begin{align*}
 & \left[a\otimes t^{m},[b\otimes t^{n},u\otimes t^{k}]\right]\\
=\  & \frac{1}{4}(n-2k)\left[a\otimes t^{m},bu\otimes t^{n+k-1}\right]\\
=\  & \frac{1}{16}(n-2k)(m-2n-2k+2)a(bu)\otimes t^{m+n+k-2},
\end{align*}
and
\begin{align*}
 & \left[u\otimes t^{k},\left[a\otimes t^{m},b\otimes t^{n}\right]\right]\\
=\  & \frac{1}{2}(m-n)\left[u\otimes t^{k},ab\otimes t^{m+n-1}\right]\\
=\  & -\frac{1}{2}(m-n)\left[ab\otimes t^{m+n-1},u\otimes t^{k}\right]\\
=\  & -\frac{1}{8}(m-n)(m+n-2k-1)(ab)u\otimes t^{m+n+k-2}.
\end{align*}
We see that
\begin{align*}
\left[a\otimes t^{m},[b\otimes t^{n},u\otimes t^{k}]\right]+\left[b\otimes t^{n},[u\otimes t^{k},a\otimes t^{m}]\right]+\left[u\otimes t^{k},[a\otimes t^{m},b\otimes t^{n}]\right]=0
\end{align*}
if and only if
\begin{eqnarray}
 &  & (n-2k)(m-2n-2k+2)a(bu)-(m-2k)(n-2m-2k+2)b(au)\nonumber \\
 &  & \ \ \ \ \ -2(m-n)(m+n-2k-1)(ab)u=0.
\end{eqnarray}
We can show that this holds for all $m,n,k\in\Z$ if and only if $a(bu)=(ab)u=b(au)$,
noticing that
\begin{align*}
(n-2k)(m-2n-2k+2)-(m-2k)(n-2m-2k+2)-2(m-n)(m+n-2k-1)=0
\end{align*}
from a straightforward calculation.

Case 3. Consider $a\in A,\ u,v\in U$. For $m,n,k\in\Z$, we have
\begin{align*}
 & \left[a\otimes t^{m},[u\otimes t^{n},v\otimes t^{k}]_{+}\right]\\
=\  & \left[a\otimes t^{m},(u\circ v)\otimes t^{n+k}+\frac{1}{2}n(n-1)\langle u,v\rangle\delta_{n+k-1,0}K\right]\\
=\  & \frac{1}{2}(m-n-k)a(u\circ v)\otimes t^{m+n+k-1}+\binom{m}{3}\langle a,u\circ v\rangle\delta_{m+n+k-2,0}K,
\end{align*}
\begin{align*}
 & \left[v\otimes t^{k},[a\otimes t^{m},u\otimes t^{n}]\right]_{+}\\
=\  & \frac{1}{4}(m-2n)\left[v\otimes t^{k},au\otimes t^{m+n-1}\right]\\
=\  & \frac{1}{4}(m-2n)\left(\left(v\circ(au)\right)\otimes t^{k+m+n-1}+\frac{1}{2}k\left(k-1\right)\langle v,au\rangle\delta_{k+m+n-2,0}K\right)\\
=\  & \frac{1}{4}(m-2n)(v\circ(au))\otimes t^{k+m+n-1}+\frac{1}{8}(m-2n)k(k-1)\langle v,au\rangle\delta_{k+m+n-2,0}K.
\end{align*}
We see that (the super Jacobi identity)
\[
\left[a\otimes t^{m},[u\otimes t^{n},v\otimes t^{k}]_{+}\right]-\left[u\otimes t^{n},[a\otimes t^{m},v\otimes t^{k}]\right]_{+}=\left[[a\otimes t^{m},u\otimes t^{n}],v\otimes t^{k}\right]_{+},
\]
which amounts to
\[
\left[a\otimes t^{m},[u\otimes t^{n},v\otimes t^{k}]_{+}\right]-\left[v\otimes t^{k},[a\otimes t^{m},u\otimes t^{n}]\right]_{+}+\left[u\otimes t^{n},[v\otimes t^{k},a\otimes t^{m}]\right]_{+}=0
\]
holds if and only if
\begin{eqnarray}
2(m-n-k)a(u\circ v)-(m-2n)(v\circ(au))-(m-2k)(u\circ(av))=0\label{condition1}
\end{eqnarray}
and
\begin{eqnarray}
\frac{1}{6}m(m-1)(m-2)\langle a,u\circ v\rangle-\frac{1}{8}(m-2n)k(k-1)\langle v,au\rangle-\frac{1}{8}(m-2k)n(n-1)\langle u,av\rangle=0\label{condition2}
\end{eqnarray}
when $m+n+k-2=0$. Indeed, it is straightforward (by taking $m=2k$
with $n\ne k$, and then taking $m=2n$ with $n\ne k$) to show that
(\ref{condition1}) holds for all $m,n,k\in\Z$ if and only if
\begin{eqnarray}
a(u\circ v)=v\circ au=u\circ av.
\end{eqnarray}
We can also prove that relation (\ref{condition2}) holds for all
$m,n,k\in\Z$ if and only if
\begin{eqnarray}
\langle v,au\rangle=\langle u,av\rangle=\frac{4}{3}\langle a,u\circ v\rangle.
\end{eqnarray}
(The ``only if'' part follows by taking $m=0$, $n=3$, $k=-1$,
and then taking $k=0$, $m=3,\ n=-1$. The proof of ``if part''
is straightforward.)

Case 4. Let $u,v,w\in U$. For $m,n,k\in\Z$, we have
\begin{align*}
\left[u\otimes t^{m},\left[v\otimes t^{n},w\otimes t^{k}\right]_{+}\right] & =\left[u\otimes t^{m},(v\circ w)\otimes t^{n+k}\right]\\
 & =-\left[(v\circ w)\otimes t^{n+k},u\otimes t^{m}\right]\\
 & =-\frac{1}{4}(n+k-2m)(v\circ w)u\otimes t^{n+k+m-1}.
\end{align*}
Using permutations we also have
\[
\left[v\otimes t^{n},\left[w\otimes t^{k},u\otimes t^{m}\right]_{+}\right]=-\frac{1}{4}(k+m-2n)(w\circ u)v\otimes t^{n+k+m-1},
\]
\[
\left[w\otimes t^{k},\left[u\otimes t^{m},v\otimes t^{n}\right]_{+}\right]=-\frac{1}{4}(m+n-2k)(u\circ v)w\otimes t^{n+k+m-1}.
\]
We see that (the super Jacobi identity)
\[
\left[u\otimes t^{m},[v\otimes t^{n},w\otimes t^{k}]_{+}\right]+\left[v\otimes t^{n},[u\otimes t^{m},w\otimes t^{k}]_{+}\right]=\left[[u\otimes t^{m},v\otimes t^{n}]_{+},w\otimes t^{k}\right]
\]
amounts to
\[
\left[u\otimes t^{m},[v\otimes t^{n},w\otimes t^{k}]_{+}\right]+\left[v\otimes t^{n},[w\otimes t^{k},u\otimes t^{m}]_{+}\right]+\left[w\otimes t^{k},[u\otimes t^{m},v\otimes t^{n}]_{+}\right]=0,
\]
which holds if and only if
\begin{eqnarray}
(n+k-2m)(v\circ w)u+(k+m-2n)(w\circ u)v+(m+n-2k)(u\circ v)w=0.
\end{eqnarray}
It is straightforward to show that the latter relation holds for all
$m,n,k\in\Z$ if and only if
\begin{eqnarray}
(u\circ v)w=(v\circ w)u=(w\circ u)v.
\end{eqnarray}
(It can be readily seen that the ``if'' part is true, whereas the
``only if'' part can be established by taking $m=0,\ n=1,\ k=2$,
and taking $m=0,\ n=1,\ k=-1$.)
\end{proof}

From now on, we assume that $A$ is a commutative associative algebra equipped with a $\C$-valued
symmetric bilinear form $\langle\cdot,\cdot\rangle_{A}$ which is
associative in the sense that
\begin{eqnarray}
\langle ab,c\rangle_{A}=\langle a,bc\rangle_{A}\ \ \ \ \text{for }a,b,c\in A.
\end{eqnarray}
Furthermore, assume that $U$ is an $A$-module equipped with a $\C$-valued
symmetric bilinear form $\langle\cdot,\cdot\rangle_{U}$ and a symmetric
bilinear map
\begin{align}
\circ:\ U\times U\rightarrow A;\ \ (u,v)\mapsto u\circ v,
\end{align}
satisfying the following conditions for $a\in A,$ $u,v\in U$:
\begin{eqnarray}
 &  & \langle au,v\rangle_{U}=\langle u,av\rangle_{U}=\frac{4}{3}\langle a,u\circ v\rangle_{A},\label{bilinear form prop}\\
 &  & a(u\circ v)=au\circ v=u\circ av.
\end{eqnarray}
In view of Theorem \ref{iff-Lie-super-algebra},
we have a Lie superalgebra $\mathcal{L}(A,U)$.

\begin{exa}\label{NS-algebra}
{\em We here show that the Neveu-Schwarz superalgebra is a special case of $\mathcal{L}(A,U)$.
Assume that $A$ has an identity written as $\frac{1}{2}\omega$
with $\langle \omega,\omega\rangle =\frac{1}{2}$. Then
\begin{align}
[\omega_{m+1},\omega_{n+1}]
&=\frac{1}{2}(m-n)(\omega\omega)_{m+n+1}+\binom{m+1}{3}\delta_{m+n,0}\langle \omega,\omega\rangle K\nonumber\\
&=(m-n)\omega_{m+n+1}+\frac{1}{12}(m^3-m)\delta_{m+n,0} K
\end{align}
for $m,n\in \mathbb{Z}$. Thus we have a copy of the Virasoro algebra with $L(m)=\omega_{m+1}$ for $m\in \mathbb{Z}$.
Furthermore, assume that $g$ is a vector in $U$ such that $g\circ g =2\omega\in A$. From (\ref{bilinear form prop}) , we have
$$\langle g,\omega g\rangle=\frac{4}{3}\langle \omega, g\circ g\rangle=\frac{8}{3}\langle \omega, \omega\rangle
=\frac{4}{3},$$
which implies $\langle g,g\rangle=\frac{1}{2}\langle g,\omega g\rangle=\frac{2}{3}$. Then
\begin{align}
&[g_{m+1},g_{n}]=(g\circ g)_{m+n+1}+\frac{1}{2}m(m+1)\langle g,g\rangle\delta_{m+n,0}K
=2\omega_{m+n+1}+\frac{1}{3}m(m+1)\delta_{m+n,0} K\nonumber\\
&\quad =2L(m+n)+\frac{1}{3}m(m+1)\delta_{m+n,0} K
\end{align}
for $m,n\in \mathbb{Z}$. Set $G_{m-\frac{1}{2}}=g_m$ for $m\in \mathbb{Z}$. We have
\begin{align}
[G_{m+\frac{1}{2}},G_{n-\frac{1}{2}}]=2L(m+n)+\frac{1}{3}m(m+1)\delta_{m+n,0} K
\end{align}
and
\begin{align*}
&[L(m),G_{n+\frac{1}{2}}]=[\omega_{m+1},g_{n+1}]
=\frac{1}{4}\left(m+1-2(n+1)\right)(\omega g)_{m+n+1}=\frac{1}{2}(m-2n-1)g_{m+n+1}\\
&\quad =\left(\frac{1}{2}m-n-\frac{1}{2}\right)G_{m+n+\frac{1}{2}}.
\end{align*}
Under these assumptions, we obtain a copy of the Neveu-Schwarz superalgebra inside $\mathcal{L}(A,U)$.}
\end{exa}

\begin{rem}
{\em Let $A$ be any commutative associative algebra
equipped with a symmetric associative bilinear form $\langle\cdot,\cdot\rangle$.
Let $U$ be an ideal, which is naturally an $A$-module. Equip $U$
with the symmetric bilinear form inherited from $A$ and define $u\circ v=\frac{3}{4}uv$
for $u,v\in U$. Then it can be readily seen that all the conditions
in Proposition \ref{iff-Lie-super-algebra} hold. On the other hand,
for a commutative associative algebra $A$, any linear functional
$f$ on $A^{2}$ (an ideal of $A$) gives rise to a symmetric associative
bilinear form $\langle\cdot,\cdot\rangle$ on $A$ by
\[
\langle a,b\rangle=f(ab)\ \ \ \ \mbox{ for }a,b\in A.
\]
In case that $A$ has an identity, we can easily see that every symmetric
associative bilinear form on $A$ can be obtained this way.} \end{rem}

\begin{exa}
{\em Let $A$ be a finite-dimensional commutative associative
algebra generated by one element. Then $A=\C[x]/(p(x))$, where $p(x)$
is a monic polynomial. Furthermore, each ideal of $\C[x]/(p(x))$
is of the form $q(x)\C[x]/(p(x))$ where $q(x)$ is a divisor of $p(x)$.
Write
\[
p(x)=(x-\alpha_{1})^{n_{1}+1}\cdots(x-\alpha_{r})^{n_{r}+1},
\]
where $\alpha_{1},\dots,\alpha_{r}$ are distinct complex numbers
and $n_{1},\dots,n_{r}$ are nonnegative integers. We have
\begin{eqnarray}
\C[x]/(p(x))\simeq\C[x]/((x-\alpha_{1})^{n_{1}+1})\oplus\cdots\oplus\C[x]/((x-\alpha_{r})^{n_{r}+1}),
\end{eqnarray}
a direct sum of algebras. Set $A_{i}=\C[x]/((x-\alpha_{i})^{n_{i}+1})$
for $1\le i\le r$. For any associative bilinear form $\langle\cdot,\cdot\rangle$
on $A$, we have $\langle A_{i},A_{j}\rangle=0$ for $1\le i\ne j\le r$.

Consider the special case $A=\C[x]/((x-\alpha)^{n+1})$ with $\alpha\in\C,\ n\ge0$.
Then $A=\oplus_{i=0}^{n}\C a^{i}$ with $a=\overline{x-\alpha}$.
Let $f$ be any linear functional on $A$ such that $f(a^{n})=1$.
Use $f$ to define a symmetric associative bilinear form $\langle\cdot,\cdot\rangle$
on $A$ by $\langle u,v\rangle=f(uv)$ for $u,v\in A$. For $0\le i,j\le n$,
we have
\[
\langle a^{i},a^{n-i}\rangle=f(a^{n})=1\ \ \mbox{ and }\ \langle a^{i},a^{j}\rangle=f(a^{i+j})=0
\]
if $i+j\ge n+1$. We see that the corresponding bilinear form on $A$
is non-degenerate. In fact, every non-degenerate symmetric associative
bilinear form up to a scalar multiple can be obtained this way. Note
that $\C[x]/(x^{n+1})\simeq\C[x]/((x-\alpha)^{n+1})$ for any $\alpha\in\C$.
Then it follows that there always exist non-degenerate symmetric associative
bilinear forms on $\C[x]/(p(x))$.}
\end{exa}

\section{Vertex operator superalgebra $V_{\mathcal{L}(A,U)}(\ell,0)$ and
its irreducible modules}


In this section, for any complex number $\ell$, we associate a vertex
operator superalgebra $V_{\mathcal{L}(A,U)}(\ell,0)$ to the Lie superalgebra $\mathcal{L}(A,U)$,
and we classify irreducible $V_{\mathcal{L}(A,U)}(\ell,0)$-modules.

Let $A$ be a commutative associative algebra with a symmetric and
associative bilinear form $\langle\cdot,\cdot\rangle$ and let $U$
be an $A$-module equipped with a symmetric bilinear form $\langle\cdot,\cdot\rangle$
and a symmetric bilinear map $U\times U\rightarrow A$, such that
the conditions in Theorem \ref{iff-Lie-super-algebra} hold. Recall
from Section 4 that we have a Lie superalgebra
\begin{eqnarray}
\mathcal{L}(A,U)=(A+U)\otimes\C[t,t^{-1}]\oplus\C K,
\end{eqnarray}
where the even part and odd part are
\[
\mathcal{L}(A,U)^{0}=A\otimes\C[t,t^{-1}]\oplus\C K,\quad\mathcal{L}(A,U)^{1}=U\otimes\C[t,t^{-1}].
\]
It can be readily seen that $\mathcal{L}(A,U)$ is a $\frac{1}{2}\mathbb{Z}$-graded
algebra with $\text{deg}K=0$ and
\begin{eqnarray}
\deg\left(a\otimes t^{n}\right)=-n+1,\ \ \ \ \deg\left(u\otimes t^{n}\right)=-n+\frac{1}{2}
\end{eqnarray}
for $a\in A,\ u\in U,\ n\in\mathbb{Z}$.

For $u\in A\oplus U$, $n\in\mathbb{Z}$, following the usual practice,
alternatively set
\[
u(n)=u\otimes t^{n}\in\mathcal{L}(A,U).
\]
For $u\in A+U$, form a generating function
\begin{eqnarray}
u(x)=\sum_{n\in\mathbb{Z}}u(n)x^{-n-1}\in\mathcal{L}(A,U)[[x,x^{-1}]].
\end{eqnarray}
Then the defining commutator relations in Proposition \ref{iff-Lie-super-algebra}
can be written as
\begin{align}
[a(x_{1}),b(x_{2})] & =(ab)\left(x_{2}\right)\frac{\partial}{\partial x_{2}}\left(x_{1}^{-1}\delta\left(\frac{x_{2}}{x_{1}}\right)\right)+\frac{1}{2}\left(\frac{\partial}{\partial x_{2}}(ab)(x_{2})\right)x_{1}^{-1}\delta\left(\frac{x_{2}}{x_{1}}\right)\nonumber \\
 & \quad+\frac{1}{6}\left(\frac{\partial}{\partial x_{2}}\right)^{3}\left(x_{1}^{-1}\delta\left(\frac{x_{2}}{x_{1}}\right)\right)\langle a,b\rangle K,\\{}
[a(x_{1}),u(x_{2})] & =\frac{3}{4}(au)(x_{2})\frac{\partial}{\partial x_{2}}x_{1}^{-1}\delta\left(\frac{x_{2}}{x_{1}}\right)+\frac{1}{2}\left(\frac{\partial}{\partial x_{2}}(au)(x_{2})\right)x_{1}^{-1}\delta\left(\frac{x_{2}}{x_{1}}\right),\\{}
[u(x_{1}),v(x_{2})]_{+} & =(u\circ v)\left(x_{2}\right)x_{1}^{-1}\delta\left(\frac{x_{2}}{x_{1}}\right)+\frac{1}{2}\left(\frac{\partial}{\partial x_{2}}\right)^{2}x_{1}^{-1}\delta\left(\frac{x_{2}}{x_{1}}\right)\langle u,v\rangle K
\end{align}
for $a,b\in A,\ u,v\in U$.

Set
\begin{eqnarray}
\mathcal{L}(A,U)^{+}=(A\oplus U)\otimes\mathbb{C}[t]\oplus\mathbb{C}K,\ \ \ \
\mathcal{L}(A,U)^{-}=(A\oplus U)\otimes t^{-1}\mathbb{C}[t^{-1}],
\end{eqnarray}
which are subalgebras with $\mathcal{L}(A,U)=\mathcal{L}(A,U)^{+}\oplus\mathcal{L}(A,U)^{-}$
as a vector space.

Let $\ell\in\mathbb{C}$. Denote by $\mathbb{C}_{\ell}$ the one-dimensional
$\mathcal{L}(A,U)^{+}$-module $\C$ with $K$ acting as scalar $\ell$
and with $(A+U)\otimes\mathbb{C}[t]$ acting trivially. (Note that
$\mathcal{L}(A,U)^{+}=(A+U)\otimes\mathbb{C}[t]\oplus\mathbb{C}K$
is a direct sum of algebras.) Form an induced module
\begin{eqnarray}
V_{\mathcal{L}(A,U)}(\ell,0)=U(\mathcal{L}(A,U))\otimes_{U\left(\mathcal{L}(A,U)^{+}\right)}\mathbb{C}_{\ell}.
\end{eqnarray}
In view of the P-B-W theorem, we have $V_{\mathcal{L}(A,U)}(\ell,0)=U(\mathcal{L}(A,U)^{-})$
as a $U(\mathcal{L}(A,U)^{-})$-module. Furthermore, we have
\begin{eqnarray}
V_{\mathcal{L}(A,U)}(\ell,0)\simeq S(L(A)^{-})\otimes\Lambda(L(U)^{-}),
\end{eqnarray}
where $L(A)^{-}=A\otimes t^{-1}\mathbb{C}[t^{-1}]$, $L(U)^-=U\otimes t^{-1}\mathbb{C}[t^{-1}]$,
$S(\cdot)$ and $\Lambda(\cdot)$ denote the symmetric algebra
and the exterior algebra, respectively.

Set
\[
\boldsymbol{1}=1\otimes1\in V_{\mathcal{L}(A,U)}(\ell,0).
\]
Identify $A\oplus U$ as a subspace of $V_{\mathcal{L}(A,U)}(\ell,0)$
through the linear map
\begin{align}
a+u\mapsto a(-1)\boldsymbol{1}+u(-1)\boldsymbol{1}\ \ \text{for\ }a\in A,\ u\in U.
\end{align}
Since $\mathcal{L}(A,U)^{+}$ is a $\frac{1}{2}\mathbb{Z}$-graded
subalgebra, $V_{\mathcal{L}(A,U)}(\ell,0)$ becomes a $\frac{1}{2}\mathbb{Z}$-graded
$\mathcal{L}(A,U)$-module by defining $\deg\boldsymbol{1}=0$. We
have
\[
V_{\mathcal{L}(A,U)}(\ell,0)_{(n)}=0\ \ \ \mbox{ for }n\in\frac{1}{2}\Z\mbox{ with }n<0,
\]
\[
V_{\mathcal{L}(A,U)}(\ell,0)_{(\frac{1}{2})}=0=V_{\mathcal{L}(A,U)}(\ell,0)_{(1)},
\]
and
\begin{eqnarray}
V_{\mathcal{L}(A,U)}(\ell,0)_{(0)}=\C{\bf 1},\ \ \ \ V_{\mathcal{L}(A,U)}(\ell,0)_{(2)}=A,\ \ \ \ V_{\mathcal{L}(A,U)}(\ell,0)_{(\frac{3}{2})}=U.
\end{eqnarray}

Define a linear operator $D$ on Lie superalgebra $\mathcal{L}(A,U)$
by
\[
D(K)=0,\ \ \ \ D(v(n))=-nv(n-1)\ \ \ \mbox{ for }v\in A+U,\ n\in\Z.
\]
In terms of generating functions we have
\[
D(v(x))=\frac{d}{dx}v(x)\ \ \ \mbox{ for }v\in A+U.
\]
It is straightforward to see that $D$ is a derivation of $\mathcal{L}(A,U)$.
Then $D$ naturally acts on $U(\mathcal{L}(A,U))$ as a derivation.
It is clear that $D$ preserves the subalgebra $\mathcal{L}(A,U)^{+}$.
Consequently, $D$ acts on $V_{\mathcal{L}(A,U)}(\ell,0)$ with $D(X\otimes1)=D(X)\otimes1$
for $X\in U(\mathcal{L}(A,U))$, such that
\[
D{\bf 1}=0,\ \ \ \ [D,v(x)]=\frac{d}{dx}v(x)\ \ \ \mbox{ for }v\in A+U.
\]
Equip $ V_{\mathcal{L}(A,U)}(\ell,0)$ with the $\mathbb{Z}_2$-graded structure with
\begin{eqnarray}
V_{\mathcal{L}(A,U)}(\ell,0)^{\bar 0}=\oplus_{n\in\Z}V_{\mathcal{L}(A,U)}(\ell,0)_{(n)},\ \ \ \ V_{\mathcal{L}(A,U)}(\ell,0)^{\bar 1}=\oplus_{n\in\frac{1}{2}+\Z}V_{\mathcal{L}(A,U)}(\ell,0)_{(n)}.
\end{eqnarray}
By the standard arguments in \cite{Li} (cf. \cite{FKRW}, \cite{MP},
\cite{KW}, \cite{Kac}) we have:

\begin{prop}
There exists a vertex superalgebra structure on $V_{\mathcal{L}(A,U)}(\ell,0)$
which is uniquely determined by the condition that $\boldsymbol{1}$
is the vacuum vector and
\begin{eqnarray}
Y(v,x)=v(x)=\sum_{n\in\mathbb{Z}}v(n)x^{-n-1}\quad\text{ for }v\in A\oplus U.
\end{eqnarray}
\end{prop}

Note that as $Y(v,x)=\sum_{n\in\Z}v_{n}x^{-n-1}$ in the standard
vertex operator notation, we have
\begin{eqnarray}
v_{n}=v(n)\ \ \ \mbox{ for }v\in A+U,\ n\in\Z.
\end{eqnarray}
Since $V_{\mathcal{L}(A,U)}(\ell,0)=U(\mathcal{L}(A,U)^{-}){\bf 1}$, we see that
vertex superalgebra $V_{\mathcal{L}(A,U)}(\ell,0)$ is strongly generated by its subspace $A+ U$.
Note that
\begin{align*}
&\deg a_n=\deg a(n)=1-n=\deg a-n-1\quad \text{ for }a\in A,\\
&\deg u_n=\deg u(n)=\frac{1}{2}-n=\deg u-n-1\quad \text{ for }u\in U.
\end{align*}
Then it follows that $V_{\mathcal{L}(A,U)}(\ell,0)$ is a $\frac{1}{2}\mathbb{Z}$-graded vertex superalgebra.

\begin{rem}
{\em Recall the linear operator $\mathcal{D}$ on $V_{\mathcal{L}(A,U)}(\ell,0)$
with $\mathcal{D}(v)=v_{-2}{\bf 1}$ for $v\in V_{\mathcal{L}(A,U)}(\ell,0)$.
For $a,b\in A=V_{\mathcal{L}(A,U)}(\ell,0)_{(2)}$, we have $a_0b=\frac{1}{2}\mathcal{D}(ab)\in \mathcal{D}(A)$ as
\begin{align}
a_0b=a_0b_{-1}{\bf 1}=b_{-1}a_{0}{\bf 1}+[a_0,b_{-1}]{\bf 1}=0+\frac{1}{2}(ab)_{-2}{\bf 1}
=\frac{1}{2}\mathcal{D}(ab),
\end{align}
where we use the bracket relation (\ref{am-bn-lie}).}
\end{rem}




Next, we classify irreducible $\frac{1}{2}\N$-graded $V_{\mathcal{L}(A,U)}(\ell,0)$-modules.
Let $W$ be an $\mathcal{L}(A,U)$-module. For $u\in A+U$, set
\[
u_{W}(x)=\sum_{n\in\mathbb{Z}}u(n)x^{-n-1}\in\left(\text{End}W\right)[[x,x^{-1}]].
\]

An $\mathcal{L}(A,U)$-module $W$ is said to be {\em restricted}
if for any $w\in W,\ u\in A+U$, $u(n)w=0$ for all sufficiently
positive integers $n$. Assume that $W$ is a {\em $\frac{1}{2}\N$-graded $\mathcal{L}(A,U)$-module}
in the sense that $W$ is an $\mathcal{L}(A,U)$-module equipped with
a $\frac{1}{2}\mathbb{N}$-grading $W=\oplus_{n\in\frac{1}{2}\mathbb{N}}W(n)$
such that
\begin{align}
\mathcal{L}(A,U)_{(p)}W(q)\subset W(p+q)\ \ \ \mbox{ for }p,q\in\frac{1}{2}\mathbb{Z}.
\end{align}
It is clear that a $\frac{1}{2}\N$-graded $\mathcal{L}(A,U)$-module
is necessarily a restricted $\mathcal{L}(A,U)$-module.

Using a method similar to that in \cite{Li} we have:

\begin{prop}\label{voa-module-Lie-module}
For any restricted $\mathcal{L}(A,U)$-module
$W$ of level $\ell\in\C$, there exists a $V_{\mathcal{L}(A,U)}(\ell,0)$-module
structure $Y_{W}(\cdot,x)$ which is uniquely determined by $Y_{W}(u,x)=u_{W}(x)$
for $u\in A+U$. On the other hand, for any $V_{\mathcal{L}(A,U)}(\ell,0)$-module
$(W,Y_{W})$, $W$ is a restricted $\mathcal{L}(A,U)$-module of level
$\ell$ with $u_{W}(x)=Y_{W}(u,x)$ for $u\in A+U$.
\end{prop}

\begin{rem}\label{classical-fact}
{\em Let $\g=\oplus_{q\in\frac{1}{2}\Z}\g(q)$
be any $\frac{1}{2}\Z$-graded Lie superalgebra with even part $\g^{0}=\oplus_{n\in\Z}\g(n)$
and odd part $\g^{1}=\oplus_{q\in\frac{1}{2}+\Z}\g(q)$. Note that
$\g(0)$ is a Lie subalgebra of $\g^{0}$. It is clear that for any
(irreducible) $\frac{1}{2}\N$-graded $\g$-module $W=\oplus_{q\in\frac{1}{2}\N}W(q)$,
$W(0)$ is an (irreducible) $\g(0)$-module. On the other hand, for
any irreducible $\g(0)$-module $U$, we have an irreducible $\frac{1}{2}\N$-graded
$\g$-module $W$ with $W(0)\simeq U$ as a $\g(0)$-module. Therefore,
irreducible $\frac{1}{2}\N$-graded $\g$-modules are in one-to-one
correspondence with irreducible $\g(0)$-modules. }
\end{rem}

Recall that $\mathcal{L}(A,U)$ is a $\frac{1}{2}\Z$-graded Lie superalgebra
with $\deg K=0$, 
\begin{eqnarray*}
\deg a(n)=1-n\ \ \mbox{ for }a\in A,\ n\in\Z,\text{ and }\ \deg u(r)=\frac{1}{2}-r\ \ \mbox{ for }u\in U,\ r\in\Z.
\end{eqnarray*}
Then $\mathcal{L}(A,U)(0)=A\otimes\C t+\C K$. For $a,b\in A,\ m\in\Z$,
from the bracket relation (\ref{am-bn-lie}) (see Theorem \ref{iff-Lie-super-algebra})
we have $[a(m),b(m)]=0$. Thus $\mathcal{L}(A,U)(0)$ is an abelian Lie algebra.
Combining Proposition \ref{voa-module-Lie-module} with Remark \ref{classical-fact} we
immediately have:

\begin{thm}\label{irred-module-1-1}
Let $S(A)$ be the symmetric algebra over $A$. Then
irreducible $\frac{1}{2}\N$-graded $V_{\mathcal{L}(A,U)}(\ell,0)$-modules
are in one-to-one correspondence with irreducible $S(A)$-modules.
\end{thm}

Note that Theorem \ref{irred-module-1-1} indicates (one can show) that the Zhu algebra of $V_{\mathcal{L}(A,U)}(\ell,0)$
is isomorphic to $S(A)$. As we work on $\C$, if $A$ is of countable dimension,
every irreducible $S(A)$-module is one-dimensional and uniquely determined by a linear functional on $A$.

Set
\begin{eqnarray}
\mathcal{L}(A,U)_{(-)}=\sum_{n\in\frac{1}{2}\Z_{+}}\mathcal{L}(A,U)_{(-n)}=A\otimes t^{2}\C[t]+U\otimes t\C[t],
\end{eqnarray}
a subalgebra of $\mathcal{L}(A,U)$. Recall that $\mathcal{L}(A,U)_{(0)}=A\otimes\C t+\C K$.
The Borel subalgebra $\mathcal{L}(A,U)_{(-)}+A\otimes \C t+\C K$ is
a direct sum of algebras $\mathcal{L}(A,U)_{(-)}+A\otimes\C t$ and
$\C K$, while $\mathcal{L}(A,U)_{(-)}$ is an ideal of $\mathcal{L}(A,U)_{(-)}+A\otimes\C t$.

Let $\ell\in\C,\ \lambda\in A^{*}$ (a linear functional). Denote
by $\C_{\ell,\lambda}$ the module for the Borel subalgebra $\mathcal{L}(A,U)_{(-)}+A\otimes \C t+\C K$,
where $\C_{\ell,\lambda}=\C$ as a vector space on which $\mathcal{L}(A,U)_{(-)}$
acts trivially, $K$ acts as scalar $\ell$, and $a(1)$ acts as scalar
$\lambda(a)$ for $a\in A$. Form an induced module (a generalized
Verma module)
\begin{eqnarray}
M(\ell,\lambda)=U(\mathcal{L}(A,U))\otimes_{U(\mathcal{L}(A,U)_{(-)}+A\otimes t+\C K)}\C_{\ell,\lambda}.
\end{eqnarray}
(Note that $V_{\mathcal{L}(A,U)}(\ell,0)$ is a quotient module of
$M(\ell,0)$.) Again, by the P-B-W theorem we have
\begin{eqnarray}
M(\ell,\lambda)\simeq S(L(A)_{(+)})\otimes\Lambda(L(U)_{(+)}),\label{eq:M(l,lambda)}
\end{eqnarray}
where $L(A)_{(+)}=A\otimes\C[t^{-1}]$ and $L(U)_{(+)}=U\otimes\C[t^{-1}]$.
Define $\deg\C_{\ell,\lambda}=0$ to make $M(\ell,\lambda)$ a $\frac{1}{2}\N$-graded
$\mathcal{L}(A,U)$-module. Since the degree-zero subspace, which
equals $\C_{\ell,\lambda}$, is one-dimensional, $M(\ell,\lambda)$
has a unique maximal submodule. Denote by $L(\ell,\lambda)$ the irreducible
quotient module modulo the unique maximal graded submodule of $M(\ell,\lambda)$.
In conclusion, we have:

\begin{coro}\label{irreducible-modules}
Up to equivalence, $L(\ell,\lambda)$
with $\lambda\in A^{*}$ exhaust all the irreducible $\frac{1}{2}\N$-graded
$V_{\mathcal{L}(A,U)}(\ell,0)$-modules.
\end{coro}

Now, we assume that $A$ has an identity written as $\frac{1}{2}\omega$ with $\langle \omega,\omega\rangle=\frac{1}{2}$.
From Example \ref{NS-algebra}, we have
\begin{align}
[L(m),L(n)]=(m-n)L(m+n)+\frac{1}{12}(m^3-m)\delta_{m+n,0}\ell
\end{align}
for $m,n\in \mathbb{Z}$, where $Y(\omega,z)=\sum_{n\in \mathbb{Z}}L(n)z^{-n-2}$.
From the bracket relations (\ref{am-bn-lie})  and (\ref{am-un-lie})  in Theorem \ref{iff-Lie-super-algebra}, we also have
\begin{align}
&[L(m),a_n]=(m+1-n)a_{m+n}+\frac{1}{6}(m^3-m)\langle \omega,a\rangle \delta_{m+n-1,0}\ell,\\
&[L(m),u_n]=\frac{1}{2}(m+1-2n)u_{m+n}
\end{align}
for $a\in A,\ u\in U$. In particular, we have
\begin{align}
&[L(-1),a_n]=-na_{n-1},\quad [L(0),a_n]=(1-n)a_n,\\
&[L(-1),u_n]=-nu_{n-1},\quad [L(0),u_n]=\frac{1}{2}(1-2n)u_n,
\end{align}
which imply
\begin{align*}
&L(0)a=L(0)a_{-1}{\bf 1}=2a_{-1}{\bf 1}=2a\  \  \text{ for }a\in A=V_{\mathcal{L}(A,U)}(\ell,0)_{(2)},\\
&L(0)u=L(0)u_{-1}{\bf 1}=\frac{3}{2} u_{-1}{\bf 1}=\frac{3}{2}u\   \  \text{ for }u\in U=V_{\mathcal{L}(A,U)}(\ell,0)_{(\frac{3}{2})}.
\end{align*}
It then follows that $\omega$ is a conformal vector of $V_{\mathcal{L}(A,U)}(\ell,0)$ with central charge $\ell$
and $V_{\mathcal{L}(A,U)}(\ell,0)$ is a conformal vertex superalgebra with
\begin{align}
V_{\mathcal{L}(A,U)}(\ell,0)_{(2)}=A,\  \   V_{\mathcal{L}(A,U)}(\ell,0)_{(\frac{3}{2})}=U
\end{align}
such that $V_{\mathcal{L}(A,U)}(\ell,0)$ is strongly generated by $A+U$.
If both $A$ and $U$ are assumed to be finite-dimensional, it follows from the P-B-W theorem that
all the homogeneous subspaces of $V_{\mathcal{L}(A,U)}(\ell,0)$ are finite-dimensional, so that $V_{\mathcal{L}(A,U)}(\ell,0)$
is a vertex operator superalgebra.

Note that from Proposition \ref{Inv-form-V} there exists a unique symmetric invariant bilinear form
$\langle\cdot,\cdot\rangle$ on $V_{\mathcal{L}(A,U)}(\ell,0)$ with $\langle {\bf 1},{\bf 1}\rangle=1$.
For $a,b\in A=V_{\mathcal{L}(A,U)}(\ell,0)_{(2)}$, we have
\begin{align}
a_3b=a_3b_{-1}{\bf 1}=b_{-1}a_3{\bf 1}+[a_3,b_{-1}]{\bf 1}
=b_{-1}a_3{\bf 1}+ 2(ab)_1{\bf 1}+\langle a,b\rangle_A K {\bf 1}=\ell \langle a,b\rangle_A {\bf 1},
\end{align}
where $\langle\cdot,\cdot\rangle_A$ denotes the given bilinear form on $A$.
Similarly, for $u,v\in U=V_{\mathcal{L}(A,U)}(\ell,0)_{(\frac{3}{2})}$, we have
\begin{align}
u_2v=u_2v_{-1}{\bf 1}=(u\circ v)_1{\bf 1}+\langle u,v\rangle_U K{\bf 1}=\ell \langle u,v\rangle_U{\bf 1}.
\end{align}
As $V_{\mathcal{L}(A,U)}(\ell,0)_{(0)}=\C {\bf 1}$,
$V_{\mathcal{L}(A,U)}(\ell,0)$ has a unique maximal (graded) ideal $J$.
If the bilinear forms $\langle\cdot,\cdot\rangle_A$ on $A$ and $\langle\cdot,\cdot\rangle_U$ on $U$
are both non-degenerate and if $\ell$ is nonzero, then
it follows from Remark \ref{bilinear-forms-A-U} that the maximal ideal $J$ which coincides with
the kernel of the invariant bilinear form on $V_{\mathcal{L}(A,U)}(\ell,0)$
intersect with $A$ and $U$ trivially. Set
\begin{align}
L_{\mathcal{L}(A,U)}(\ell,0)=V_{\mathcal{L}(A,U)}(\ell,0)/J,
\end{align}
a simple conformal vertex superalgebra.  Thus we have proved:

\begin{thm}\label{existence-simple-vosa}
Let $\ell$ be a nonzero complex number and assume that the given bilinear forms on $A$ and $U$ are both non-degenerate.
Then $L_{\mathcal{L}(A,U)}(\ell,0)$ is a simple conformal vertex superalgebra with
\begin{align}
L_{\mathcal{L}(A,U)}(\ell,0)_{(2)}=A,\quad L_{\mathcal{L}(A,U)}(\ell,0)_{(\frac{3}{2})}=U,
\end{align}
such that $L_{\mathcal{L}(A,U)}(\ell,0)$ is strongly generated by $A+U$.
Furthermore, if both $A$ and $U$ are finite-dimensional, then $L_{\mathcal{L}(A,U)}(\ell,0)$ is a vertex operator superalgebra
(with finite-dimensional homogeneous subspaces).
\end{thm}

 Recall that for $\lambda\in A^{*}$, $L(\ell,\lambda)$ is an irreducible $\frac{1}{2}\N$-graded (weak) $V_{\mathcal{L}(A,U)}(\ell,0)$-module
 with $L(\ell,\lambda)(0)=\C_{\ell,\lambda}\ (=\C)$. Set $v_{\lambda}=1\in \C_{\ell, \lambda}$.
 Note that $L(0)v_{\lambda}=\lambda(\omega)v_{\lambda}$.
 Then it follows that
 \begin{align}
 L(\ell,\lambda)=\oplus_{q\in \frac{1}{2}\N}L(\ell,\lambda)_{(\lambda(\omega)+q)},
 \end{align}
 where $L(\ell,\lambda)_{(\mu)}=\{ v\in L(\ell,\lambda)\ |\ L(0)v=\mu v\}$ for $\mu\in \C$.
 Note that if both $A$ and $U$ are finite-dimensional, then it follows from the P-B-W theorem that
 all the homogeneous subspaces are finite-dimensional, namely, $L(\ell,\lambda)_{(\mu)}$ is an (ordinary)
 $V_{\mathcal{L}(A,U)}(\ell,0)$-module.

 Follow \cite{FHL} to define the restricted dual
 \begin{align}
  L(\ell,\lambda)'=\oplus_{q\in \frac{1}{2}\N}(L(\ell,\lambda)_{(\lambda(\omega)+q)})^{*}.
\end{align}
From \cite{Xu}, $L(\ell,\lambda)'$ is also a $\frac{1}{2}\N$-graded (weak) $V_{\mathcal{L}(A,U)}(\ell,0)$-module,
where the vertex operator map
 $$Y'(\cdot,z):\  V\rightarrow ({\rm End} L(\ell,\lambda)')[[z,z^{-1}]]$$
 is given by
\begin{align}
\langle Y'(v,z)\alpha,w\rangle=\left\langle \alpha,Y\left( e^{zL(1)}z^{-2L(0)}e^{\pi L(0)}v,z^{-1}\right)w\right\rangle
\end{align}
for $v\in V,\ \alpha\in L(\ell,\lambda)',\ w\in L(\ell,\lambda)$. For $a\in A=V_{\mathcal{L}(A,U)}(\ell,0)_{(2)}$,
from (\ref{invariance-an}) we have
\begin{align}
\langle a_n\alpha,w\rangle=\langle \alpha,a_{2-n}w\rangle\quad \text{ for }n\in \Z.
\end{align}
Let $v_{\lambda}^*\in (\C_{\ell,\lambda})^*\subset L(\ell,\lambda)'$ defined by $\langle v_{\lambda}^*,v_{\lambda}\rangle =1$.
For $a\in A$, we have
$$\langle a_1 v_{\lambda}^{*},v_{\lambda}\rangle=\langle v_{\lambda}^{*},a_1v_{\lambda}\rangle
=\lambda(a)\langle v_{\lambda}^{*},v_{\lambda}\rangle=\lambda(a),$$
which implies $a_1v_{\lambda}^{*}=\lambda(a)v_{\lambda}^{*}$. In conclusion, we have:

\begin{prop}\label{form-on-module}
Assume that both $A$ and $U$ are finite-dimensional. Then for every $\lambda\in A^{*}$,
the contragredient dual $L(\ell,\lambda)'$ of $V_{\mathcal{L}(A,U)}(\ell,0)$-module $L(\ell,\lambda)$ is isomorphic to
$L(\ell,\lambda)$ itself and there is a non-degenerate symmetric invariant bilinear form on $L(\ell,\lambda)$.
\end{prop}

\vskip10pt

\noindent\footnotesize{\textbf{H. Li}:  Department of Mathematical Sciences, Rutgers University, Camden, NJ 08102 USA; \texttt{hli@camden.rutgers.edu}}

\noindent\footnotesize{\textbf{N. Yu}: School of Mathematical
Sciences, Xiamen University, Fujian, 361005, CHINA; \texttt{
ninayu@xmu.edu.cn}}


\begin{thebibliography}{FKRW}

\bibitem[BLP]{BLP}
C. Bai, H.-S. Li, Y. Pei, $\phi_{\epsilon}$-coordinated modules for vertex algebras,
{\em J. Algebra} \textbf{426} (2015), 211--242.

\bibitem[Bo]{Bor}
R. Borcherds, Vertex algebras, Kac-Moody algebras,
and the Monster, {\em Proc. Nat. Acad. Sci. U.S.A.} \textbf{83}
(1986), no. 10, 3068--3071.

\bibitem[DLM]{DLM}
C. Dong, H.-S. Li, G. Mason, Regularity of rational
vertex operator algebras, {\em Adv. Math.} (1997).

\bibitem[FFR]{FFR}
A. J. Feingold, I. B. Frenkel, J. F. X. Ries,
Spinor construction of vertex operator algebras, triality, and $E_{8}^{(1)}$,
{\em Contemp. Math.} \textbf{121}, Amer. Math. Soc., Providence,
RI, 1991.

\bibitem[FKRW]{FKRW}
E. Frenkel, V. Kac, A. Radul, and W. Wang, $W_{1+\infty}$
and $W(\mathfrak{gl}_{N})$ with central charge $N$, {\em Commun.
Math. Phys.} \textbf{170} (1995), 337--357.

\bibitem[FHL]{FHL}
I. B. Frenkel, Y.-Z. Huang, J. Lepowsky, Axiomatic
Approach to Vertex Operator Algebras, {\em Amer. Math. Soc.}, Providence,
RI, 1993; Preprint, 1989.

\bibitem[FZ]{FZ}
I. B. Frenkel, Y.-C. Zhu, Vertex operator algebras
associated to representations of affine and Virasoro algebras, {\em
Duke Math. J.} \textbf{66} (1992), no. 1, 123--168.

\bibitem[K]{Kac}
V. Kac,  Vertex Algebras for Beginners,
Second edition, University Lecture Series 10, American Mathematical
Society, Providence, RI, 1998.

\bibitem[KW]{KW}
V. Kac, W. Wang, Vertex operator superalgebras and
their representations, Mathematical aspects of conformal and topological
field theories and quantum groups (South Hadley, MA, 1992), 161--191,
{\em Contemp. Math.}, 175, Amer. Math. Soc., Providence, RI, 1994.

\bibitem[KL]{KL99}
M. Karel and H.-S. Li, Certain generting subspaces for vertex operator algebras,
{\em J. Algebra} \textbf{217} (1999), 393--421.

\bibitem[La1]{Lam-thesis}
C.H. Lam, On the structure of vertex operator algebras and their weight two subspaces,
Ph.D. thesis, The Ohio State University, 1996.

\bibitem[La2]{Lam}
C.H. Lam, Construction of vertex operator algebras
from commutative associative algebras, {\em Commun. Algebra} \textbf{24}(14)  (1996),
4339--4360.

\bibitem[Li1]{Li-thesis}
H.-S. Li, Representation theory and tensor
product theory of vertex operator algebras, Ph.D. dissertation, Rutgers
University, New Brunswick, 1994.

\bibitem[Li2]{Li-form}
H.-S. Li, Symmetric invariant bilinear forms on vertex operator algebras,
{\em J. Pure Applied Algebra} {\bf 96} (1994), 279--297.

\bibitem[Li3]{Li}
H.-S. Li, Local systems of vertex operators, vertex
superalgebras and modules, {\em J. Pure Applied Algebra} \textbf{109}
(1996), 143--195.

\bibitem[LY]{LY}
H.-S. Li, G. Yamaskuna, On certain vertex algebras
and their modules associated with vertex algebroids, {\em J. Algebra}
\textbf{283}:1 (2005), 367--398.

\bibitem[Lia]{Lian}
B.-H. Lian, On the classification of simple vertex operator algebras,
{\em Commun. Math. Phys.} {\bf 163} (1994), 307--357.

\bibitem[MP]{MP}
A. Muerman, M. Primc, Vertex operator algebras and
representations of affine Lie algebras, {\em Acta Appl. Math.} \textbf{ 44}
(1-2), 207--215.

\bibitem[P]{Pr}
M. Primc, Vertex algebras generated by Lie algebras, {\em J. Pure Applied Algebra}
{\bf 135} (1999) 253--293.

\bibitem[S]{S} N. Scheithauer, Vertex algebras, Lie algebras, and
superstrings.\emph{ J. Algebra} \textbf{200} (1998), no. 2, 363--403.

\bibitem[X1]{Xu}
X. Xu,  Introduction to Vertex Operator Superalgebras
and Their Modules, Mathematics and its applications, Kluwer Academic
Publishers, 1998.

\bibitem[X2]{Xu-00}
X. Xu, Quadratic conformal superalgebras, {\em J. Algebra} {\bf 231} (2000) 1-38.

\bibitem[Z1]{Zhu}
Y. Zhu, Vertex operator algebras, elliptic functions
and modular forms, Ph.D. dissertation, Yale University, 1990.

\bibitem[Z2]{Zhu2}
Y. Zhu, Modular invariance of characters of
vertex operator algebras, {\em J. Amer. Math. Soc.} \textbf{9}
(1996), 237-302.
\end{thebibliography}
\end{document}